\documentclass[final,3p,twocolumn]{elsarticle} 

\usepackage{float}
\usepackage{amsmath,amssymb,amsfonts,amstext,amsmath}
\usepackage{epsfig}
\usepackage{setspace}
\usepackage{hhline}
\usepackage{color,soul}
\usepackage[dvipsnames]{xcolor}
\definecolor{Light}{gray}{.90}
\usepackage[justification=centering]{caption}
\usepackage{multirow}
\usepackage{booktabs}
\usepackage{bm}
\usepackage[english]{babel}
\usepackage{tikz}
\usepackage{tikzscale}
\usepackage{epstopdf}
\usepackage{adjustbox}
\usepackage{longtable}
\usepackage{accents}
\usepackage{lineno}
\usepackage{multirow}
\usepackage{colortbl}
\usepackage[linesnumbered,ruled,vlined]{algorithm2e}
\usepackage{ragged2e}
\usepackage{amsmath,amssymb}
\makeatletter
\usepackage{breqn}
\usepackage{amsthm}
\usepackage{pdflscape}
\usepackage{afterpage}
\usepackage{capt-of}

\newcolumntype{x}[1]{>{\raggedright\arraybackslash\hspace{0pt}}p{#1}}
\newtheorem{problem}{Problem}
\newtheorem{assumption}{Assumption}

\newtheorem{lemma}{Lemma}
\newtheorem{theorem}{Theorem}
\newtheorem{remark}{Remark}

\newcommand{\iter}[1]{\ensuremath{\langle #1 \rangle}}
\newcommand{\itrpt}[3]{\ensuremath{\bm{#1}^{\iter{#2}}_{#3}}}
\newcommand{\itrval}[3]{\ensuremath{#1^{\iter{#2}}_{#3}}}

\newcommand{\idxpt}[3]{\ensuremath{\bm{#1}^{(#2)}_{#3}}}

\SetCommentSty{mycommfont}
\SetKwInput{KwInput}{Input}        
\SetKwInput{KwOutput}{Output}       

\DeclareMathOperator*{\argmin}{arg\,min}
\DeclareMathOperator*{\argmax}{arg\,max}

\begin{document}

\begin{frontmatter}
\title{SMGO: A Set Membership Approach\\to Data-Driven Global Optimization}

\author{Lorenzo Sabug Jr.\corref{corr}}\ead{lorenzojr.sabug@polimi.it}
\cortext[corr]{Corresponding author. The first author would like to acknowledge the support of the Department of Science and Technology--Science Education Institute (DOST-SEI) of the Philippines for his research.}
\author{Fredy Ruiz}\ead{fredy.ruiz@polimi.it}
\author{Lorenzo Fagiano\corref{ack2}}\ead{lorenzo.fagiano@polimi.it}
\cortext[ack2]{This research has been supported by the Italian Ministry of University and Research (MIUR) under the PRIN 2017 grant n. 201732RS94 "Systems of Tethered Multicopters".}
\address{Dipartimento di Elettronica, Informazione e Bioingegneria, Politecnico di Milano\\ Piazza Leonardo da Vinci 32, 20133 Milano, Italy}

\begin{abstract}

Many science and engineering applications feature non-convex optimization problems where the objective function can not be handled analytically, i.e. it is a black box. Examples include design optimization via experiments, or via costly finite elements simulations. To solve these problems, global optimization routines are used. These iterative techniques must trade-off exploitation close to the current best point with exploration of unseen regions of the search space. In this respect, a new global optimization strategy based on a Set Membership (SM) framework is proposed. Assuming Lipschitz continuity of the cost function, the approach employs SM concepts to decide whether to switch from an exploitation mode to an exploration one, and vice-versa. The resulting algorithm, named SMGO (Set Membership Global Optimization) is presented. Theoretical properties regarding convergence and computational complexity are derived, and implementation aspects are discussed. Finally, the SMGO performance is evaluated on a set of benchmark non-convex problems and compared with those of other global optimization approaches.

\end{abstract}

\end{frontmatter}

\section{Introduction}
In many science and engineering fields, such as mechanical design, fluid-dynamics, electromagnetics, multi-physics simulations, control systems tuning, and chemical experiments, black-box optimization problems arise. Black-box functions are named as such due to the fact that an explicit mathematical model is unavailable, or too complicated to be handled analytically. In these cases, the optimization strategy can rely only on function values obtained through empirical tests. Moreover, these problems may present several local minima, and the time and resources required to carry out a single function evaluation are rather large, so that the solution method shall make the most efficient use of the available trials. To highlight this fact, we refer to the process of performing a new test as \textit{long function} evaluation.

Seeking a minimizer in this framework is referred to as black-box optimization~\cite{Jones1998}, or derivative-free optimization~\cite{Pham2011}. This is in contrast for example to gradient-based and Newton-type methods where the (first, or even further) derivatives are also assumed to be available. If the black-box function is assumed differentiable, one could resort to gradient estimation techniques to still employ these methods, which however require a rather large number of long function evaluations to estimate a local quantity (i.e., the gradient at the currently evaluated point)~\cite{Amaran2016} and are designed to converge to a local optimum by always improving from the initial guess, not to explore the decision space searching for a global solution. 

Due to time/resource limitations, in global optimization there are two conflicting aspects that must be considered when choosing the next test point for a long function evaluation~\cite{Amaran2016}. The first, \textit{exploitation}, pertains to improving the result by testing more points in the vicinity of the current best known one. The second, \textit{exploration}, aims to gain more information about the cost function in other regions of the search space, seeking other, possibly better, minima. Exploitation and exploration strategies are at the core of most global optimization approaches, including stochastic global search, response surface, and Lipschitz-based methods.

Stochastic global search algorithms, such as random search methods~\cite{GAO2018},  particle swarm algorithms~\cite{Kennedy2014}, and its offshoot variants~\cite{Liu2017,Han2017}, use a population of search points per iteration. According to the corresponding function values, heuristics are applied to calculate the search point locations for the next iteration (or ``generation''). These methods are widespread but require numerous long function evaluations per generation, which can be impractical. Moreover, aspects such as guaranteed convergence and optimality gap may be difficult to analyze. Response surface methods, on the other hand, generate an approximation of the black-box function via, e.g., Gaussian process regression~\cite{Jones1998}, radial basis functions~\cite{Powell2002,Bemporad2020}, or neural networks~\cite{Wang2019}. These approximations can then be used for the exploitation or exploration routines, whichever is appropriate. Albeit the number of long function evaluation is lower in this case, optimality gap is again not easily assessed, and there is an additional exponential computational burden to derive and refine the approximated objective function. 

Another category of black-box optimizers are so-called Lipschitz-based algorithms. They rely on the assumption that the long function is Lipschitz continuous. An early method, proposed independently by Piyavskii~\cite{Piyavskii1972} and Shubert~\cite{Shubert1972}, uses a known Lipschitz constant to draw the lower bounds of the black-box function given the existing samples, starting with the bounds (corners) of the search space. In this context, the Piyavskii-Shubert (``sawtooth'') method chooses the next point by searching the location of the minimum lower bound. However, the method is proposed for one-dimensional function optimization only; furthermore, Lipschitz constants are rarely known in practical cases, and need to be estimated in the process.

A later proposal named DIRECT~\cite{Jones1993, Finkel2006, KHONG2013} is a modification to Piyavskii-Shubert method, removing the need for a known Lipschitz constant and starting samples from the corners. Instead, it invariably performs a first sample at the center of the search space, which is afterwards subdivided into three intervals. The Lipschitz constant estimate is used to draw the lower bounds, which are used to select potentially optimal intervals (hyperrectangles). Such hyperrectangles are iteratively subdivided, and their corresponding centers are sampled by batch, and potential intervals for sampling are selected again. This cycle proceeds until the maximum number of iterations is reached, with each iteration performing multiple long function evaluations. An offshoot method named DISIMPL~\cite{Paulavicius2014a, Paulavicius2014b} employs simplices instead of hyperrectangles, thus extending the usage from rectangular search spaces only to more general polytopic ones. On one hand, DIRECT and DISIMPL are completely deterministic and reproducible (assuming no noise on function evaluation), however, they only admit a fixed starting point, and do not propose how to handle existing sampled points at the beginning of the algorithms. Furthermore, the batch-based sampling means that the next set of trial points is not chosen until the current batch of sampling is finished, which can mean overlooking possible precious information added by individual samples as they come in.

LIPO and AdaLIPO~\cite{Malherbe2017} use a random binary variable to decide between exploitation and exploration; the only difference being that the former assumes a fixed Lipschitz constant, while the latter estimates it from existing data. For exploitation, they choose a random point within the set of potential optimizing points, i.e. those whose lower bound on the black-box function is below the current best sample. For exploration, it chooses a random point from the entire search space. These algorithms, owing to the simple implementation, are fast in terms of choosing the next sampling point, and can easily scale to high dimensional problems without significant hit on computational burden. However, such a completely randomized exploitation/exploration sampling may be inefficient in finding new local minima, because it does not fully utilize the shape of lower bounds. Furthermore, the results of such an algorithm are unrepeatable because of its randomized design.

In the described context, we propose a new Lipschitz-based algorithm that exploits Set Membership (SM) nonlinear function approximation theory~\cite{MILANESE20112141}. So far, SM theory proved effective for nonlinear system identification~\cite{Milanese2004}, filter design~\cite{6220238}, and controller design~\cite{Fagiano2016}, among others. Here, we employ SM theory for the first time in global optimization. Our paper delivers the following contributions:

\begin{itemize}
  \item The Set Membership Global Optimization (SMGO) algorithm for Lipschitz-continuous black-box functions is introduced. It employs the concepts of SM lower- and upper bounds in exploitation and exploration routines. Unlike the above-mentioned Lipschitz-based algorithms, we address the exploration problem systematically, by casting it as minimization of the uncertainty as inferred using SM techniques. Furthermore, our proposed algorithm is completely repeatable, and allows one to consider custom starting points and initial samples when available;
  \item Theoretical properties of SMGO are derived, including convergence, optimality gap, and computational complexity;
  \item Practical aspects are discussed, and ways to improve the computational efficiency are described. In particular, an iterative implementation is described, which exploits results from previous iterations to alleviate the computations at the current one;
  \item The SMGO algorithm is compared with other representative optimizers with different test functions, showing that the performance of SMGO is very competitive w.r.t. the state of the art, especially on functions with many global minima.
\end{itemize}

A preliminary version of this work appeared in \cite{SaRF2020}. With respect to that contribution, the algorithm presented here is more efficient, moreover the theoretical results and implementation analysis are new, and more benchmark results are presented. The MATLAB code for SMGO, as well as the other scripts used to generate the results, are available on GitHub under the following URL: \verb|https://github.com/lorenzosabugjr/smgo|.

This paper is organized in seven sections. Section~\ref{sec:prob-state} gives the general problem statement and assumptions. Section~\ref{sec:smgo-algo} describes the proposed SMGO algorithm. The convergence properties and optimality gap calculations are discussed in Section~\ref{sec:algo-analysis}, with implementation notes in Section~\ref{sec:implement-notes}. Section~\ref{sec:perf-test} compares the performance of the proposed algorithm with other global optimization techniques on representative test functions, and Section~\ref{sec:conclusion} concludes this paper and provides insights on future directions.

\section{Problem Statement}
\label{sec:prob-state}

Consider a cost function $z=f_o(\bm{x}),\;f_o:\mathcal{X}\rightarrow\mathbb{R}$, where $x\in\mathcal{X}$ is the vector of decision variables,  $\mathcal{X}\subset\mathbb{R}^D$ is a compact and convex polytope (``search set''). No analytical form of $f_o$ is assumed available. The \emph{a priori} knowledge about $f_o$ is given by the following assumption: 

\begin{assumption}
\label{ass:lipschitz}
$f_o$ is Lipschitz continuous with unknown Lipschitz constant $\gamma_o$, i.e.,
\[ f_o \in \mathcal{F}(\gamma_o)\]
where

\begin{multline*}
    \mathcal{F}(\gamma_o) \doteq \Big\{ f \in \bm{C}^0(\mathcal{X}): |f(\bm{x}_1) - f(\bm{x}_2)| ~\leq \\ \gamma_o\|\bm{x}_1 - \bm{x}_2\|, \forall\bm{x}_1,\bm{x}_2 \in \mathcal{X} \Big\}
\end{multline*}

\end{assumption}

Being a Lipschitz continuous function on a compact set, $f_o$ presents a global minimum $z^*$:
\begin{equation}
  z^* = \min_{\bm{x} \in \mathcal{X}} f_o(\bm{x}).
  \label{eqn:global-minimum}
\end{equation}
Let us denote with
\begin{equation}
  \mathcal{X}^* = \left\{ ~\bm{x}\in\mathcal{X}~|~f_o(\bm{x})=z^*~ \right\}
  \label{eqn:global-minimzer}
\end{equation}
the corresponding set of global minimizers.

We further assume that it is possible to acquire information about $f_o$ by sampling (long function evaluation):

\begin{assumption}\label{ass:function-eval}
Given a vector of decision variables  $\bm{x}\in\mathcal{X}$, it is possible to sample the cost function without noise:
 \[z = f_o(\bm{x}).\]
\end{assumption}
$\,$\\
The set of collected data is denoted as:
\begin{equation}\label{eqn:data-set}
\bm{X}^{\iter{n}}= \left\{(\bm{x}^{(1)},z^{(1)}); (\bm{x}^{(2)},z^{(2)}); \ldots; (\bm{x}^{(n)},z^{(n)}) \right\}
\end{equation} 
where $n\in\mathbb{N}$ is the number of data points and $z^{(i)}=f_o(\bm{x}^{(i)})$. For simplicity, with a slight abuse of notation, we will also write $\idxpt{x}{i}{} \in \bm{X}^{\iter{n}}$ when the pair $(\bm{x}^{(i)},z^{(i)})$ belongs to the data set \eqref{eqn:data-set}.

We denote with $(\bm{x}^{*\iter{n}},z^{*\iter{n}})$ the best pair in $\bm{X}^{\iter{n}}$, where

\begin{equation}\label{eqn:best-point}
(\bm{x}^{*\iter{n}},z^{*\iter{n}})=\argmin_{(\bm{x}^{(k)},z^{(k)})\in\bm{X}^{\iter{n}}}z^{(k)}
\end{equation}

\noindent If the result of \eqref{eqn:best-point} is not unique, a lexicographic criterion is used to sort the minimizers, and the first one is picked.
We further denote with $\delta^{\iter{n}}$ the optimality gap:

\begin{equation}
\delta^{\iter{n}}=z^{*\iter{n}}-z^*.
\end{equation}

We assume that a starting number $n_0\geq1$ of data points is also available, for example collected by the user in a first testing campaign, forming the set $\bm{X}^{\iter{n_0}}\subset\mathcal{X}$.
$\,$\\

We can now state the problem addressed in this paper.
$\,$\\
\begin{problem}
\label{prob:globalopt}
Design an algorithm that, under Assumptions \ref{ass:lipschitz}-\ref{ass:function-eval}, generates a sequence of points $\{\bm{x}^{(n_0+1)},\bm{x}^{(n_0+2)},\ldots\},\,\bm{x}^{(i)}\in\mathcal{X},$ such that:
\begin{equation}
  \forall \epsilon>0,\;\exists n_\epsilon<\infty: z^{*\iter{n_\epsilon}} \leq z^* + \epsilon.
  \label{eqn:prob-state}
\end{equation}
Moreover, for a finite sequence length $N$, provide a method to compute a bound on the optimality gap $\delta^{\iter{N}}$.
\end{problem}

\section{Set Membership Global Optimization (SMGO): Algorithm}
\label{sec:smgo-algo}

SMGO addresses Problem \ref{prob:globalopt} with a sequential procedure, as common in global optimization~\cite{Jones1998,Jones1993,Malherbe2017}. At each iteration $n\geq n_0$, the next test point $\bm{x}^{(n+1)}$ is chosen with a strategy that exploits the prior knowledge on function $f_o$, given by Assumption \ref{ass:lipschitz}, and the current data set $\bm{X}^{\iter{n}}$. In particular, a valid (i.e., consistent with data) estimate of the Lipschitz constant $\gamma_o$ is derived at $n_0$ and updated at each $n> n_0$, in order to estimate upper and lower bounds on $f_o$. SMGO then leverages the latter to compute $\bm{x}^{(n+1)}$, choosing between an exploitation mode and an exploration one.

\subsection{Preliminaries: Lipschitz constant estimation, cost function bounds, and search mesh}
\label{subsec:mode-sampling-update}
At the start of each iteration ($n> n_0$), one long function evaluation $z^{(n)} = f_o(\bm{x}^{(n)})$ is carried out, where $\bm{x}^{(n)}$ is the test point selected at the previous iteration. This resulting new data point ($\bm{x}^{(n)},~z^{(n)}$) is added to the set of collected samples:
\[
\bm{X}^{\iter{n}}=\bm{X}^{\iter{n-1}}\cup (\bm{x}^{(n)},~z^{(n)}). 
\]
At the first iteration, $n=n_0$, an initial estimate of the Lipschitz constant, $\gamma^{\iter{n_0}}$, is computed as:

\begin{equation}\label{eqn:gamma-starting}
\gamma^{\iter{n_0}} = \max_{i,j=1,\ldots,n_0,\,i\neq j} \frac{| z^{(i)} - z^{(j)} |}{ \| \bm{x}^{(i)} - \bm{x}^{(j)} \|}
\end{equation}

If $n_0=1$, the Lipschitz constant estimate can be initialized as an arbitrarily small number.
When $n>n_0$, the minimum feasible Lipschitz constant estimate $\gamma^{\iter{n}}$ is updated as follows:

\begin{equation}
  \gamma^{\iter{n}} = \max \left(\gamma^{\iter{n-1}}, \max_{k=1,\ldots,n-1} \frac{| z^{(n)} - z^{(k)} |}{ \| \bm{x}^{(n)} - \bm{x}^{(k)} \|}\right).
  \label{eqn:gamma-update}
\end{equation}
Such an estimate is guaranteed to be compatible with the available data and prior information on $f_o$, given Assumption \ref{ass:function-eval} \cite{Milanese2004}. Moreover, $\gamma^{\iter{n}}$ is monotonically increasing with $n$ and convergence to the true Lipschitz constant $\gamma_o$ is guaranteed if the test points densely cover the search set \cite{Fagiano2016}. This is the case for the SMGO algorithm as shown in Section \ref{sec:algo-analysis}. 

Given the data set $\bm{X}^{\iter{n}}$ and Lipschitz constant estimate $\gamma^{\iter{n}}$ we can compute the following lower and upper bounds on $f_o$ ~\cite{Milanese2004}:

\begin{equation}
  \label{eqn:lower-bounds}
  \underline{z}^{\iter{n}}(\bm{x}) \triangleq \max_{k =1, \ldots, n} \left(z^{(k)} - \mu\gamma^{\iter{n}} \|\bm{x}-\bm{x}^{(k)} \|\right),
\end{equation}

\begin{equation}
  \label{eqn:upper-bounds}
  \overline{z}^{\iter{n}}(\bm{x}) \triangleq \min_{k =1 \ldots n} \left(z^{(k)} + \mu\gamma^{\iter{n}} \|\bm{x}-\bm{x}^{(k)} \|\right).
\end{equation}

\noindent When $\mu=1$, the functions $\underline{z}^{\iter{n}},\,\overline{z}^{\iter{n}}$ represent, respectively, the minimum and maximum values that the objective function can take for a given vector of decision variables $\bf{x}\in\mathcal{X}$, given the information collected up to iteration $n$.\\

In SMGO, an (user-defined) overestimation parameter $\mu > 1$ is used in \eqref{eqn:lower-bounds}-\eqref{eqn:upper-bounds} to compensate for the fact that the $\gamma^{\iter{n}}$ estimated from the currently available data points is an underestimator of the true Lipschitz constant $\gamma_o$.

\begin{figure}[!t]
	\centering
	\includegraphics[width=\columnwidth]{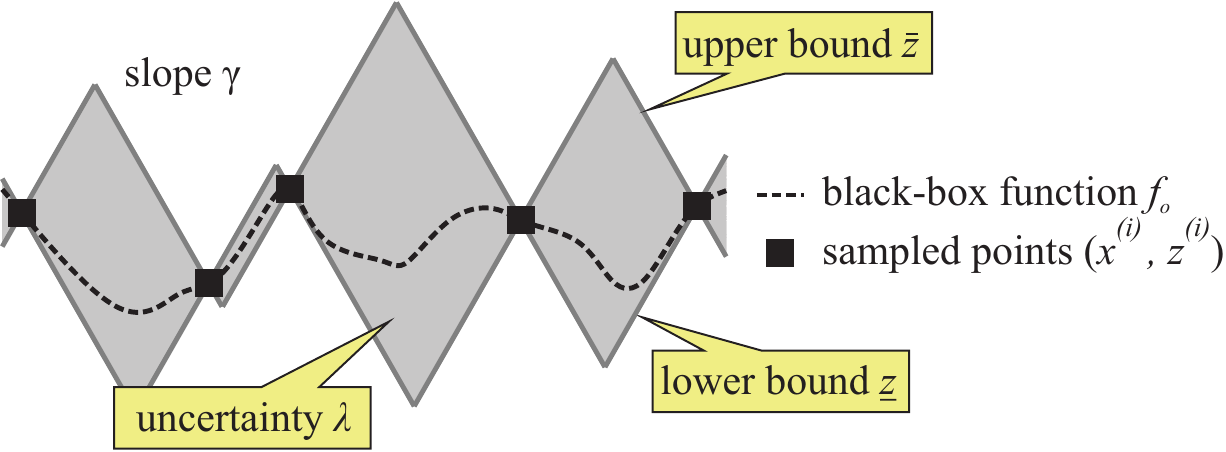}
	\caption{Scalar example of optimal lower and upper bounds $\underline{z}^{\iter{n}}(\bm{x}),\,\overline{z}^{\iter{n}}(\bm{x})$ from finite samples of $f_o$, and of the corresponding uncertainty interval $\lambda^{\iter{n}}(\bm{x})$.}
	\label{fig:sm-overview}
\end{figure}

The set 

\[
    \mathbb{F}^{\iter{n}} = \left\{ f\in\mathcal{F}(\gamma^{\iter{n}}):  \underline{z}^{\iter{n}}(\bm{x}) \leq f(\bm{x}) \leq \overline{z}^{\iter{n}}(\bm{x}) \right\}
\]

\noindent is referred to as the \textit{feasible functions set}~\cite{Milanese2004} at iteration $n$. The functions $\underline{z}^{\iter{n}}(\bm{x})$ and $\overline{z}^{\iter{n}}(\bm{x})$ are intersections of hypercones (referred simply as cones in the remainder), each one with vertex at a tested point. A representation of these bounds is depicted in Fig.~\ref{fig:sm-overview} for a function $f_o$ with scalar input argument. It should be emphasized that $\underline{z}^{\iter{n}}$ and $\overline{z}^{\iter{n}}$ are calculated from the data in $\bm{X}^{\iter{n}}$, without any need for additional evaluations of $f_o$.

The uncertainty about the value of $f_o$ at the point $\bm{x}\in\mathcal{X}$ is computed as:

\begin{equation}
  \label{eqn:sm-uncertainty}
  \lambda^{\iter{n}}(\bm{x}) = \overline{z}^{\iter{n}}(\bm{x}) - \underline{z}^{\iter{n}}(\bm{x}).
\end{equation}

The SMGO algorithm proceeds in two steps. First, an exploitation routine searches for a local minimizer of $\underline{z}^{\iter{n}}$. Then, an exploration step seeks a maximizer of $\lambda^{\iter{n}}$. Exact solutions for such subproblems can be computed by utilizing Hyperbolic Voronoi Diagrams (HVDs)~\cite{Milanese2004,Milanese2007}, however at rather high computational cost. To increase computational efficiency, we thus introduce a mesh given by segments $l_{i,j}$ connecting any two points $\bm{x}^{(i)},\,\bm{x}^{(j)} \in \bm{X}^{\iter{n}}$:

\begin{multline}
\label{eqn:segment}
    l_{i,j} = \Big\{\bm{x}=a\idxpt{x}{i}{} + (1-a)\idxpt{x}{j}{}: a \in [0, 1], \\ \idxpt{x}{i}{},\idxpt{x}{j}{} \in \bm{X}^{\iter{n}},\,i\neq j\Big\}
\end{multline}

Such a mesh allows us to derive analytic solutions to the computation of the next test point $\bm{x}^{(n+1)}$, according to the two strategies presented next.

\subsection{Exploitation (Mode~$\theta$)}
\label{subsec:mode-theta}

Consider the best data point $(\bm{x}^{*\iter{n}},z^{*\iter{n}})$ \eqref{eqn:best-point} at iteration $n$. Mode~$\theta$ searches a point $\bm{x}^{\iter{n}}_{\theta *}$ with the highest predicted improvement w.r.t. $\bm{x}^{*\iter{n}}$, according to the lower bound $\underline{z}^{\iter{n}}$ \eqref{eqn:lower-bounds}. Specifically, the exploitation routine considers the set of segments $l_{n^*,i}$ between $\bm{x}^{*\iter{n}}$ and any other sampled point $\bm{x}^{(i)}\in\bm{X}^{\iter{n}}$ to search for the candidate point $\bm{x}^{\iter{n}}_{\theta *}$. Here, index $n^*$ is the data-point index of the pair $(\bm{x}^{*\iter{n}},z^{*\iter{n}})$ in the data set $\bm{X}^{\iter{n}}$.

We first evaluate the lower bounds on $f_o(\bm{x}):\bm{x} \in l_{n^*,i}$ considering only the hypercones defined by the data points $(\bm{x}^{*\iter{n}},z^{*\iter{n}})$ and $(\bm{x}^{(i)},z^{(i)})$, i.e.:

\begin{equation}\label{eqn:lowerbound-1}
z^{*\iter{n}} - \mu\gamma^{\iter{n}} \|\bm{x}-\bm{x}^{*\iter{n}}\|
\end{equation}
and
\begin{equation}\label{eqn:lowerbound-2}
z^{(i)} - \mu\gamma^{\iter{n}} \|\bm{x}-\bm{x}^{(i)}\|.
\end{equation}

\begin{figure}[!t]
	\centering
	\includegraphics[width=\columnwidth]{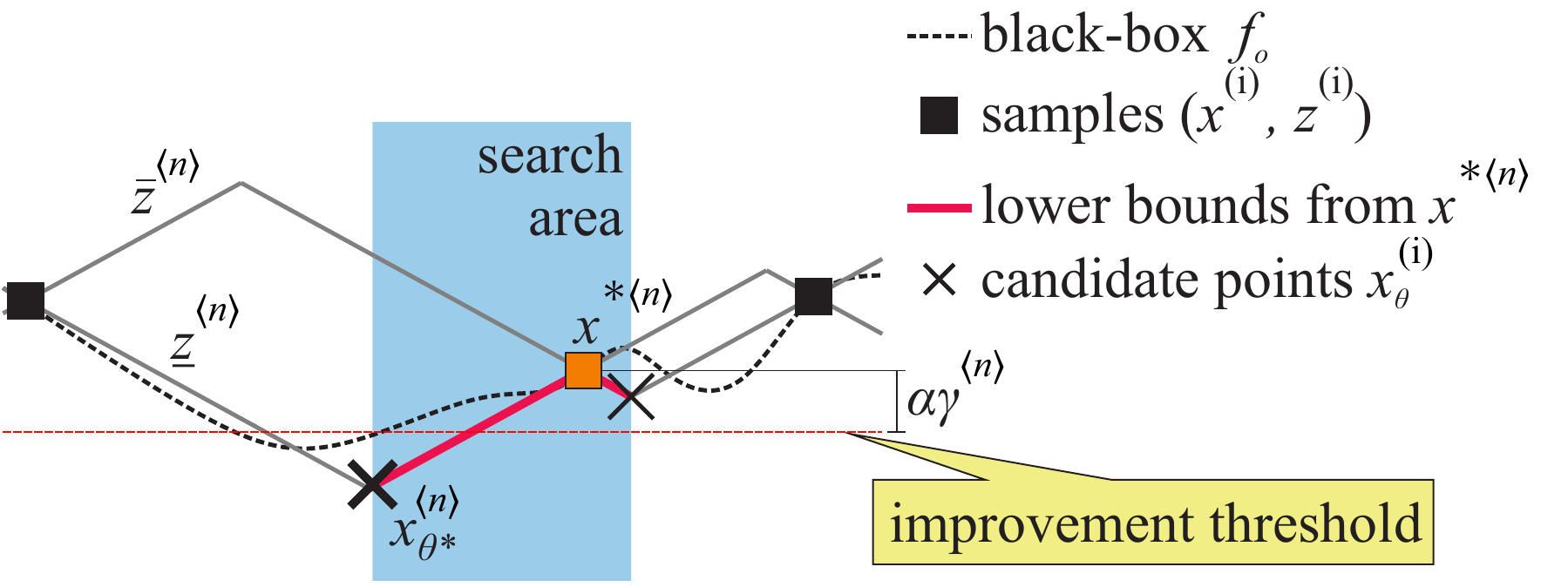}
	\caption{Qualitative example of the exploitation routine for a scalar problem.}
	\label{fig:exploit}
\end{figure}

The intersection of these two cones on the segment $l_{n^*,i}$ is the smallest feasible value of  $f_o(\bm{x}):\bm{x} \in l_{n^*,i}$ considering the information given by the two points $(\bm{x}^{*\iter{n}},z^{*\iter{n}}),\,(\bm{x}^{(i)},z^{(i)})$ and the Lipschitz constant estimate $\gamma^{\iter{n}}$. Such an intersection, denoted with $ \bm{x}_\theta^{(i)}$, can be derived analytically on the basis of geometrical considerations (see Fig.~\ref{fig:exploit}). In fact, denoting 
\begin{equation}
s_i= \frac{z^{(i)} - z^{*\iter{n}}}{\| \bm{x}^{(i)} - \bm{x}^{*\iter{n}} \|} \geq 0,
\label{eqn:exploit1-const}
\end{equation}

\noindent it follows that

\begin{equation}
  \bm{x}_\theta^{(i)} = \bm{x}^{*\iter{n}} + \frac{1 - \frac{s_i}{\mu\gamma^{\iter{n}}}}{2} (\bm{x}^{(i)} - \bm{x}^{*\iter{n}}).
  \label{eqn:exploit1-loc}
\end{equation}

By evaluating \eqref{eqn:exploit1-loc} for all possible segments, we obtain the set $\bm{X}^{\iter{n}}_\theta = \{ \bm{x}_\theta^{(i)},i=1,\ldots,n,\,i\neq n^*\}$. Then $\bm{x}^{\iter{n}}_{\theta *}$ is chosen as:

\begin{subequations}  \label{eqn:exploit-choice}
\begin{gather}
  \bm{x}^{\iter{n}}_{\theta *} = \argmin_{\bm{x} \in \bm{X}^{\iter{n}}_\theta} \underline{z}^{\iter{n}}(\bm{x})   \label{eqn:exploit-choice-a}\\
  \mathrm{s.t.}\nonumber\\
  \underline{z}^{\iter{n}}(\bm{x}) = z^{*\iter{n}} - \mu\gamma^{\iter{n}} \| \bm{x} - \bm{x}^{*\iter{n}} \|   \label{eqn:exploit-choice-b}
  \end{gather}
\end{subequations}

\noindent where the constraint \eqref{eqn:exploit-choice-b} ensures that the point selected is such that the lower bound $\underline{z}^{\iter{n}}$ considering the information provided by all the sampled points coincides with the lower bound provided by the current best point $(\bm{x}^{*\iter{n}},z^{*\iter{n}})$ alone. This usually happens in the vicinity of $\bm{x}^{*\iter{n}}$, in particular at the border of the hyperbolic Voronoi cell pertaining to $\bm{x}^{*\iter{n}}$, see \cite{Milanese2004}. The feasibility of problem \eqref{eqn:exploit-choice} is always guaranteed, as shown in Lemma \ref{lemma:exploit-feasible} in Section \ref{sec:algo-analysis}.

To decide whether $\bm{x}^{\iter{n}}_{\theta *}$ shall be eventually taken as the next test point $\bm{x}^{(n+1)}$ for the long function, or the exploration routine shall be called instead, the following condition is evaluated:

\begin{equation}
  \underline{z}^{\iter{n}}(\bm{x}^{\iter{n}}_{\theta *}) \leq z^{*\iter{n}} - \alpha \gamma^{\iter{n}},
  \label{eqn:exploit1-cond}
\end{equation}

\noindent where $\alpha \gamma^{\iter{n}}$ 
is referred to as the \textit{expected improvement threshold} and $\alpha \in [0,1)$ is the SMGO \textit{exploitation parameter}, adjusted by the user. If \eqref{eqn:exploit1-cond} is fulfilled (as in the example of Fig.~\ref{fig:exploit}), it means that the potential improvement provided by $\bm{x}^{\iter{n}}_{\theta *}$ is larger than the threshold. In this case, we set 

\[
  \bm{x}^{(n+1)} = \bm{x}^{\iter{n}}_{\theta *}.
\]

\noindent Otherwise, SMGO switches to exploration mode, described in the following subsection.

\subsection{Exploration (Mode~$\psi$)}
 
In this mode, the candidate point $\bm{x}^{\iter{n}}_{\psi *} \in \mathcal{X}$ is chosen as the one with largest $\lambda^{\iter{n}}$ value among selected candidates. Sampling $\bm{x}^{\iter{n}}_{\psi *}$ will  yield $\lambda^{\iter{n+1}}(\bm{x}^{\iter{n}}_{\psi *}) = 0$ (due to noiseless function evaluation, see Assumption \ref{ass:function-eval}), and the uncertainty in its vicinity will correspondingly decrease.

Similar to Mode~$\theta$, the search for $\bm{x}^{\iter{n}}_{\psi *}$ considers the segments among the samples in $\bm{X}^{\iter{n}}$, in order to limit the computational complexity. In this case, we consider as candidates all the midpoints between any pair $(\bm{x}^{(i)},\bm{x}^{(j)}) \in \bm{X}^{\iter{n}}$. The reason is that midpoints generally feature larger uncertainty than other points on the segments, which are closer to one of the two extremes. Hence, we form the set of midpoints:

\begin{equation}\label{eqn:explore-set}
  \bm{X}^{\iter{n}}_\psi = \left\{ \bm{x}_\psi^{(i,j)} = \frac{\bm{x}^{(i)} + \bm{x}^{(j)}}{2} ~\Big|~ \bm{x}^{(i)}, \bm{x}^{(j)} \in \bm{X}^{\iter{n}} \right\},
\end{equation}

\noindent and, according to the described rationale, the exploration routine picks the next test point as:

\begin{equation}
\label{eqn:explore-point}
  \bm{x}^{(n+1)} = \bm{x}^{\iter{n}}_{\psi *}=\argmax_{\bm{x} \in \bm{X}^{\iter{n}}_{\psi}} \lambda^{\iter{n}}(\bm{x}).
\end{equation}
Also in this case, if \eqref{eqn:explore-point} does not have a unique solution, a lexicographic criterion is added to select the maximizer.

\subsection{Algorithm summary}

A summarized flow of the SMGO method is given as pseudo-code in Algorithm~\ref{algo:smgo}.

\begin{algorithm}
  \small 
  \DontPrintSemicolon
  \KwInput{Long function $f_o$, initial data points $\bm{X}^{\iter{n_0}}$, maximum number of long function evaluations $N$, exploitation parameter $\alpha$, overestimation parameter $\mu$}
  \While{ $n_0<n\leq N$}
  {
    \tcp{Long function evaluation, data update}
    Evaluate the long function $f_o$ at $\bm{x}^{(n)}$ and measure output $z^{(n)}$. Update the data set $\bm{X}^{\iter{n}}\leftarrow \bm{X}^{\iter{n-1}}\cup (\bm{x}^{(n)}, z^{(n)})$\;
    Update current best sample $(\bm{x}^{*\iter{n}}, z^{*\iter{n}})$ and Lipschitz constant $\gamma^{\iter{n}}$ from $\bm{X}^{\iter{n}}$\;
    \tcp{Exploitation (Mode $\theta$)}
Evaluate lower bounds on segments around the current best sample $\bm{x}^{*\iter{n}}$ to choose the exploitation point $\bm{x}^{\iter{n}}_{\theta *}$\;
    \If{expected improvement condition \eqref{eqn:exploit1-cond} is met}
    {
      Assign test point for next iteration $\bm{x}^{(n+1)} \leftarrow \bm{x}^{\iter{n}}_{\theta *}$\;
    }
    \Else
    {
      \tcp{Exploration (Mode $\psi$)}
      Find $\bm{x}^{\iter{n}}_{\psi *}$ as the midpoint, over all possible segments, with largest uncertainty  $\lambda^{\iter{n}}$\;
      Assign test point for next iteration $\bm{x}^{(n+1)} \leftarrow \bm{x}^{\iter{n}}_{\psi *}$\;
    }
    Go to next iteration $n \leftarrow n+1$
  }
  Final optimal point and value: get the best sample ($\bm{x}^{*\iter{N}}, z^{*\iter{N}}$) from the set $\bm{X}^{\iter{N}}$
\caption{SMGO Algorithm}
\label{algo:smgo}
\end{algorithm}

\subsection{On the choice of $\alpha$ and $\mu$}
The SMGO parameter $\alpha$ affects only Mode~$\theta$: a smaller $\alpha$ value leads to higher tendency for exploitation. Regarding $\mu$, using a value close to one is recommended, because a higher Lipschitz constant overestimator leads to $\idxpt{x}{i}{\theta} \rightarrow \dfrac{\bm{x}^{*\iter{n}} + \idxpt{x}{i}{}}{2}$, which would not fully utilize the Set Membership approach in exploitation. At the same time, excessively large values of $\mu$ increase the conservativeness of $\lambda^{\iter{n}}$. Moreover, in the computation of $\lambda^{\iter{n}}(\bm{x}_\psi^{(i,j)}),\,\bm{x}_\psi^{(i,j)}\in\bm{X}^{\iter{n}}_\psi$ \eqref{eqn:explore-set}-\eqref{eqn:explore-point}, with larger $\mu$ values the terms $\|\bm{x}_\psi^{(i,j)}-\bm{x}^{(i)}\|,\,i,j=1,\ldots,n$
gain a higher relative importance with respect to the corresponding values of $z^{(i)}$ (see \eqref{eqn:lower-bounds}-\eqref{eqn:upper-bounds}), so that the exploration mode will tend to select the midpoints that are farthest away from the available data-points. 

\section{Theoretical Properties}
\label{sec:algo-analysis}
We now analyze the convergence and optimality properties of the SMGO algorithm. 
First, we prove that problem \eqref{eqn:exploit-choice} is always feasible, thus showing that the exploitation routine always returns a valid candidate next point $\bm{x}^{\iter{n}}_{\theta *}$.

\begin{lemma}
\label{lemma:exploit-feasible}
At any iteration $n$, problem \eqref{eqn:exploit-choice} admits at least one feasible point.
\end{lemma}

\begin{proof}
Take the pair $(\bm{x}^{(i)},z^{(i)})$ such that

\begin{equation}\label{eqn:proof-exploit-feasible-1}
\begin{array}{c}
z^{(i)}-\mu\gamma^{\iter{n}}\|\bm{x}^{*\iter{n}}-\bm{x}^{(i)}\|=\\ \max\limits_{(\bm{x},z)\in\bm{X}^{\iter{n}}\setminus(\bm{x}^{*\iter{n}},z^{*\iter{n}})} 
z-\mu\gamma^{\iter{n}}\|\bm{x}^{*\iter{n}}-\bm{x}\|,
\end{array}
\end{equation}
and consider the corresponding candidate exploitation point $\bm{x}_\theta^{(i)}$ given by  \eqref{eqn:exploit1-loc}. By construction we have that 
\[
z^{(i)}-\mu\gamma^{\iter{n}}\|\bm{x}_\theta^{(i)}-\bm{x}^{(i)}\| ~= z^{*\iter{n}}-\mu\gamma^{\iter{n}}\|\bm{x}_\theta^{(i)}-\bm{x}^{*\iter{n}}\|,
\]
and we want now to prove that this point also coincides with $\underline{z}^{\iter{n}}(\idxpt{x}{i}{\theta})$, i.e. that:

\begin{multline}
\label{eqn:proof-exploit-feasible-2}
    z^{(i)}-\mu\gamma^{\iter{n}}\|\bm{x}_\theta^{(i)}-\bm{x}^{(i)}\| ~=~ \\ \max\limits_{(\bm{x},z)\in\bm{X}^{\iter{n}}} 
    z-\mu\gamma^{\iter{n}}\|\bm{x}_\theta^{(i)}-\bm{x}\|.
\end{multline}

Assume now, for the purpose of contradiction, that:

\begin{multline}
\label{eqn:proof-exploit-feasible-3}
    \exists (\bm{x}^{(j)},z^{(j)})\neq (\bm{x}^{(i)},z^{(i)}): \\ z^{(j)}-\mu\gamma^{\iter{n}}\|\bm{x}_\theta^{(i)}-\bm{x}^{(j)}\| ~> z^{(i)}-\mu\gamma^{\iter{n}}\|\bm{x}_\theta^{(i)}-\bm{x}^{(i)}\|,
\end{multline}

\noindent thus invalidating \eqref{eqn:proof-exploit-feasible-2}. Then, denoting with $a=\frac{1}{2}(1 - \frac{s_i}{\mu\gamma^{\iter{n}}})$ (see \eqref{eqn:exploit1-loc}), we would have:

{
\small
\begin{align*}
    z^{(j)}-\mu\gamma^{\iter{n}}\|\bm{x}^{*\iter{n}}-\bm{x}^{(j)}\| &\geq\\
    z^{(j)}-\mu\gamma^{\iter{n}}\|\bm{x}^{*\iter{n}}-\bm{x}_\theta^{(i)}\|-\mu\gamma^{\iter{n}}\|\bm{x}_\theta^{(i)}-\bm{x}^{(j)}\| &>\\
    z^{(i)}-\mu\gamma^{\iter{n}}\|\bm{x}_\theta^{(i)}-\bm{x}^{(i)}\|-\mu\gamma^{\iter{n}}\|\bm{x}^{*\iter{n}}-\bm{x}_\theta^{(i)}\| &=\\
    z^{(i)}-\mu\gamma^{\iter{n}}\left(\|\bm{x}_\theta^{(i)}-\bm{x}^{(i)}\|+\|\bm{x}_\theta^{(i)}-\bm{x}^{*\iter{n}}\|\right) &=\\
    z^{(i)}-\mu\gamma^{\iter{n}}\left((1-a)\|\bm{x}^{*\iter{n}}-\bm{x}^{(i)}\|+a\|\bm{x}^{*\iter{n}}-\bm{x}^{(i)}\|\right) &=\\
    z^{(i)}-\mu\gamma^{\iter{n}}\|\bm{x}^{*\iter{n}}-\bm{x}^{(i)}\| &,
\end{align*}
}

which contradicts equation \eqref{eqn:proof-exploit-feasible-1}. Then, \eqref{eqn:proof-exploit-feasible-3} is false, while \eqref{eqn:proof-exploit-feasible-2} holds true, implying that point $\bm{x}_\theta^{(i)}$ satisfies constraint \eqref{eqn:exploit-choice-b}, hence proving the result.
\end{proof}

\subsection{Convergence of SMGO algorithm}

We next show that after finite iterations, an $\epsilon$-suboptimal point is obtained, as required by Problem \ref{prob:globalopt}. For simplicity, we consider the following technical assumption, which in practice can be replaced by a vertex estimation mechanism presented in Section \ref{subsec:coverage}. Let us denote with $V\in\mathbb{N}$ the number of vertices of the polytope $\mathcal{X}$, and with $\idxpt{c}{v}{}$ the $v$-th vertex.

\begin{assumption}
\label{assumption:assume-corners-sampled}
The starting data set $\bm{X}^{\iter{n_0}}$ contains all the vertices of $\mathcal{X}$: $\idxpt{c}{v}{} \in\bm{X}^{\iter{n_0}},\, \forall v \in [1,\ldots,V]$.
\end{assumption}

We start by proving four lemmas that are instrumental to show SMGO convergence.

In the following, for a point $\bm{x}$ and positive scalar $r$, denote with $\mathcal{B}(\bm{x},r)$ the closed hyper-ball:
\[
\mathcal{B}(\bm{x},r)=\left\{\bm{y}:\|\bm{y}-\bm{x}\|\leq r\right\}
\]

\begin{lemma}
\label{lemma:exploit-will-fail}
Mode $\theta$ will fail after a finite number of iterations.
\end{lemma}

\begin{proof} 
At iteration $n$, denote with $\mathcal{G}^{\iter{n}} $ the subset of the search space  where a candidate test point is accepted to perform an exploitation, i.e.,

\[
  \mathcal{G}^{\iter{n}} \triangleq \left\{ \bm{x} \in \mathcal{X} ~:~ \underline{z}^{\iter{n}}(\bm{x}) < z^{*\iter{n}} - \alpha\gamma^{\iter{n}} \right\}.
\]

When a candidate point inside this set is chosen as $\bm{x}^{(n+1)}$ in Mode~$\theta$, two cases arise:
\begin{enumerate}
    \item The new sample does not change the current best point $\bm{x}^{*\iter{n}}$.  In this  case,  $z^{(n+1)} \geq z^{*\iter{n}}$. Denote as $\mathcal{B}(\bm{x}^{(n+1)},r_\theta)$, the hyper-ball centered at $\bm{x}^{(n+1)}$ with finite radius 

\[
  r_\theta \triangleq \frac{z^{(n+1)} - (z^{*\iter{n}} - \alpha \gamma^{\iter{n}})}{\mu\gamma^{\iter{n+1}}}.
\]

\noindent It follows that $\forall \bm{x} \in \mathcal{B}(\bm{x}^{(n+1)},r_\theta)$ 
\[\underline{z}^{\iter{n+1}}(\bm{x}) \geq z^{*\iter{n}} - \alpha \gamma^{\iter{n}}.\]

Therefore all the points inside the hyper-ball are not eligible for exploitation. Hence, we have:

\[
  \mathcal{G}^{\iter{n+1}} = \mathcal{G}^{\iter{n}} \setminus \mathcal{B}(\bm{x}^{(n+1)},r_\theta)
\]

\noindent i.e., $\mathcal{G}$ diminishes by a finite amount, which is valid even when $\gamma^{\iter{n}}$ updates, because $\gamma^{\iter{n}}$ is bounded by $\gamma_o$.

\item  The new sample replaces the current best point $\bm{x}^{*\iter{n}}$. In this case $z^{(n+1)} < z^{*\iter{n}}$ and therefore $\mathcal{G}^{\iter{n+1}} \subset \mathcal{G}^{\iter{n}}$ due to the new threshold. Moreover, the hyper-ball $\mathcal{B}(\bm{x}^{(n+1)},r_\theta)$, with $r_\theta = \dfrac{\alpha}{\mu}$ around the new sample is also removed from $\mathcal{G}$. 
\end{enumerate}
Hence, $\mathcal{G}^{\iter{n}} = \varnothing$ after finite iterations, and Mode~$\theta$ will fail, proving the lemma.
\end{proof}

\begin{lemma}
\label{lemma:line-directed-towards}
    Let Assumption \ref{assumption:assume-corners-sampled} hold and consider any unsampled point $\hat{\bm{x}}$ in the interior of $\mathcal{X}$ and any sample $\idxpt{x}{i}{} \in \bm{X}^{\iter{n}}$. Consider the open half-space:
    \[
    \mathcal{Q} = \Big\{ \bm{x} \in \mathcal{X} ~:~ (\bm{x} - \idxpt{x}{i}{})^\top (\hat{\bm{x}} - \idxpt{x}{i}{}) > 0 \Big\}.
    \]
    Then, there exists at least one sample $\idxpt{x}{j}{} \in \bm{X}^{\iter{n}}\cap\mathcal{Q},\,i\neq j$.
\end{lemma}

\begin{proof}
    \noindent Let us assume, for the purpose of contradiction, that the lemma does not hold. That is, there is no sample in the open half-space $\mathcal{Q}$. Due to Assumption \ref{assumption:assume-corners-sampled}, this would also apply to all vertices $\idxpt{c}{v}{},\,v=1,\ldots,V$ of $\mathcal{X}$. Therefore, we would have:
    \begin{align}
    \label{eqn:dot-product-leq-zero}
        (\hat{\bm{x}} - \idxpt{x}{i}{})^\top (\idxpt{c}{v}{} - \idxpt{x}{i}{}) &\leq 0,\,\forall v=1,\ldots,V.
    \end{align}
    
    \noindent Two cases could then arise:
    
    \begin{enumerate}
        \item Inequality in \eqref{eqn:dot-product-leq-zero} applies strictly for at least one vertex $\idxpt{c}{v}{}$. In this case, any such vertex $\idxpt{c}{v}{}$ and $\hat{\bm{x}}$ would reside on opposite sides of the hyper-plane
    
        \[
            \mathcal{P} \triangleq \left\{ \bm{x} ~:~ (\hat{\bm{x}} - \idxpt{x}{i}{})^\top (\bm{x} - \idxpt{x}{i}{}) = 0 \right\}.
        \]
    
        However, consider any $\bm{w} \in \mathcal{X}$. Since the latter is a convex polytope, we can write
        
        \begin{equation}
        \label{eqn:linear-comb}
            \bm{w} = \sum_{v=1}^{V} b_v \idxpt{c}{v}{}
        \end{equation}
        
        \noindent with $0\leq b_v\leq 1,\,v=1,\ldots,V$, and $\sum\limits_{v=1}^{V} b_v = 1$. Hence, it would apply that
        
        \begin{align*}
            \forall \bm{w}\in\mathcal{X},\;(\hat{\bm{x}} - \idxpt{x}{i}{})^\top(\bm{w} - \idxpt{x}{i}{}) &= \\
            (\hat{\bm{x}} - \idxpt{x}{i}{})^\top\left(\sum_{v=1}^{V} b_v \idxpt{c}{v}{} - \sum_{v=1}^{V} b_v\idxpt{x}{i}{}\right) &= \\
            \sum_{v=1}^{V} b_v (\hat{\bm{x}} - \idxpt{x}{i}{})^\top(\idxpt{c}{v}{} - \idxpt{x}{i}{}) &< 0,
        \end{align*}
        
        \noindent which would mean that all $\bm{w} \in \mathcal{X}$ reside on the opposite side of $\mathcal{P}$ with respect to $\hat{\bm{x}}$ . This can only happen if $\hat{\bm{x}} \notin \mathcal{X}$, falling in a contradiction with the assumptions.
    
        \item Equality in \eqref{eqn:dot-product-leq-zero} applies for all $\idxpt{c}{v}{},\,v=1,\ldots, V$. This would mean that all vertices (and consequently, all $\bm{w}\in\mathcal{X}$) belong to the hyper-plane $\mathcal{P}$ containing $\idxpt{x}{i}{}$. In turn, this would imply that either $\hat{\bm{x}} \notin \mathcal{X}$ (if $\|\hat{\bm{x}}-\idxpt{x}{i}{}\| > 0$), or $\hat{\bm{x}} = \idxpt{x}{i}{}$, both of which fall in contradiction with the assumptions.
    \end{enumerate} 
    Since both cases 1. and 2. lead to contradiction, the claim of the lemma must be true, thus completing the proof.
\end{proof}

\begin{lemma}
\label{lemma:theoretical-min-lambda}
    Consider any unsampled point $\hat{\bm{x}}$ and any radius $r$ such that $\bm{X}^{\iter{n}}\cap\mathcal{B}(\hat{\bm{x}},r)=\varnothing$. Then, given Assumption~\ref{ass:function-eval} and a Lipschitz constant estimate $\gamma^{\iter{n}}$ obtained with \eqref{eqn:gamma-starting}-\eqref{eqn:gamma-update}, it applies that
    \begin{equation}
    \label{eqn:min-lambda}
        \lambda^{\iter{n}}(\hat{\bm{x}}) \geq 2(\mu-1)\gamma^{\iter{n}} r.
    \end{equation}
\end{lemma}

\begin{proof}
    Consider the point $\hat{\bm{x}}$ and denote with $(\bm{x}^{(a)},z^{(a)}),\,(\bm{x}^{(b)},z^{(b)})\in\bm{X}^{\iter{n}}$ the data pairs that determine the upper and lower bounds $\overline{z}^{\iter{n}}(\hat{\bm{x}})$, $\underline{z}^{\iter{n}}(\hat{\bm{x}})$, respectively, i.e:
    \begin{multline*}
        \overline{z}^{\iter{n}}(\hat{\bm{x}}) = \min\limits_{k =1 \ldots n} \left(z^{(k)} + \mu\gamma^{\iter{n}} \|\hat{\bm{x}}-\bm{x}^{(k)} \|\right) = \\ z^{(a)} + \mu\gamma^{\iter{n}} \|\hat{\bm{x}}-\bm{x}^{(a)} \|
    \end{multline*}
    \begin{multline*}
        \underline{z}^{\iter{n}}(\hat{\bm{x}}) = \max\limits_{k =1, \ldots, n} \left(z^{(k)} - \mu\gamma^{\iter{n}} \|\hat{\bm{x}}-\bm{x}^{(k)} \|\right) = \\ z^{(b)} - \mu\gamma^{\iter{n}} \|\hat{\bm{x}}-\bm{x}^{(b)} \|.
    \end{multline*}
Then, considering that $\bm{X}^{\iter{n}}\cap\mathcal{B}(\hat{\bm{x}},r)=\varnothing$, that $(\mu-1)>0$, and taking into account \eqref{eqn:gamma-starting}-\eqref{eqn:gamma-update} we have:

{\small
\[
\begin{array}{rcl}
    \lambda^{\iter{n}}(\hat{\bm{x}}) &=& \overline{z}^{\iter{n}}(\hat{\bm{x}})-\underline{z}^{\iter{n}}(\hat{\bm{x}})\\
    &=& z^{(a)} - z^{(b)} +  \\
    &~& \mu\gamma^{\iter{n}}\left(\|\hat{\bm{x}}-\bm{x}^{(a)}\|+\|\hat{\bm{x}}-\bm{x}^{(b)} \|\right)\\
    &\geq& -\gamma^{\iter{n}}\|\bm{x}^{(a)}-\bm{x}^{(b)}\| + \\
    &~& \mu\gamma^{\iter{n}} \left(\|\hat{\bm{x}}-\bm{x}^{(a)}\|+\|\hat{\bm{x}}-\bm{x}^{(b)} \|\right)\\
    &\geq& -\gamma^{\iter{n}}\left(\|\hat{\bm{x}}-\bm{x}^{(a)}\|+\|\hat{\bm{x}}-\bm{x}^{(b)} \|\right)+ \\
    &~& \mu\gamma^{\iter{n}} \left(\|\hat{\bm{x}}-\bm{x}^{(a)}\|+\|\hat{\bm{x}}-\bm{x}^{(b)} \|\right)\\
    &\geq& 2\left(\mu-1\right)\gamma^{\iter{n}} r
\end{array}
\]
}
which proves the result.
\end{proof}

Note that in Lemma \ref{lemma:theoretical-min-lambda} we are not assuming that the estimate $\gamma^{\iter{n}}$ is larger than $\gamma_o$: the result follows by how the estimate is computed, combined with the use of $\mu>1$.

\begin{lemma}
\label{lemma:dense-pts}
    Let Assumption \ref{assumption:assume-corners-sampled} hold, and assume that Mode~$\psi$ is undertaken infinitely often as $n\rightarrow+\infty$. Then, for any $\hat{\bm{x}}\in\mathcal{X}$ and any $\sigma > 0$, $\exists n_\sigma < \infty$ such that
    
    \[ 
        \min_{\bm{x}^{(i)} \in \bm{X}^{\iter{n_\sigma}}} \| \bm{x}^{(i)} - \hat{\bm{x}} \| < \sigma. 
    \]
\end{lemma}

\begin{proof}
Consider any point $\hat{\bm{x}} \in \mathcal{X}$. If $\hat{\bm{x}}\in\bm{X}^{\iter{n}}$ for some $n<\infty$, we have $\min\limits_{\bm{x}^{(i)} \in \bm{X}^{\iter{n}}} \| \bm{x}^{(i)} - \hat{\bm{x}} \|=0$ and the claim is trivially proven. Consider then the case $\hat{\bm{x}}\notin\bm{X}^{\iter{n}},\,\forall n\in[1,+\infty)$.  Pick the nearest sample 

\begin{equation}\label{eqn:lemma-nearest-sample}
    \bar{\bm{x}}^{\iter{n}} = \min_{\bm{x}^{(i)} \in \bm{X}^{\iter{n}}} \| \bm{x}^{(i)} - \hat{\bm{x}} \|
\end{equation}

\noindent and consider the open half-space:

\[
    \mathcal{Q} = \Big\{ \bm{x} \in \mathcal{X} ~:~ (\bm{x} - \bar{\bm{x}}^{\iter{n}})^\top (\hat{\bm{x}} - \bar{\bm{x}}^{\iter{n}}) > 0 \Big\}.
\]

\noindent Now, we select sample $\tilde{\bm{x}}^{\iter{n}}$ 

\begin{equation}\label{eqn:lemma-nearest-sample2}
    \tilde{\bm{x}}^{\iter{n}} = \argmin_{\idxpt{x}{j}{} \in \bm{X}^{\iter{n}} \cap \mathcal{Q}} \| \bm{x}^{(j)} - \bar{\bm{x}}^{\iter{n}} \|.
\end{equation}
\noindent Such a sample is always guaranteed to exist, in force of Lemma \ref{lemma:line-directed-towards}. We note that due to our selection criterion for $\tilde{\bm{x}}^{\iter{n}}$, the set

\[
    \mathcal{H} = \Big\{ \bm{x} \in \mathcal{Q} ~:~ \| \bm{x} - \bar{\bm{x}}^{\iter{n}} \| < \| \tilde{\bm{x}}^{\iter{n}} - \bar{\bm{x}}^{\iter{n}} \| \Big\}
\]

\noindent does not have any sample in its interior, otherwise any such sample should have been selected as $\tilde{\bm{x}}^{\iter{n}}$ in \eqref{eqn:lemma-nearest-sample2}. The geometry of $\mathcal{H}$ is an open half-hyper-ball, with center at $\bar{\bm{x}}^{\iter{n}}$ and radius $\| \tilde{\bm{x}}^{\iter{n}} - \bar{\bm{x}}^{\iter{n}} \|$.

We now show that a midpoint inside $\mathcal{H}$ is going to be sampled after finite iterations. There exists at least the midpoint $\hat{\bm{m}} \triangleq \dfrac{\bar{\bm{x}}^{\iter{n}} + \tilde{\bm{x}}^{\iter{n}}}{2}$ in the interior of $\mathcal{H}$. Consider the ball $\mathcal{B}(\hat{\bm{m}},h)$ centered at $\hat{\bm{m}}$ with radius $h$:

\begin{equation}\label{eqn:midpoint-h}
\begin{array}{c}
    h = \max\limits_{a\in\mathbb{R}^+} a\\
    \text{s.t.}\\
    \mathcal{B}(\hat{\bm{m}},a) \subset \mathcal{H}
\end{array}
\end{equation}

\noindent Then, recalling that there are no samples inside $\mathcal{H}$, by Lemma~\ref{lemma:theoretical-min-lambda} we have:

\begin{equation}\label{eqn:midpoint-proof}
    \lambda^{\iter{n}}(\hat{\bm{m}}) \geq 2(\mu-1)\gamma^{\iter{n}} h.
\end{equation}
Furthermore, we consider the set of points not belonging to $\mathcal{H}$, whose uncertainty  is larger than that of $\hat{\bm{m}}$
\[
    \mathcal{R}^{\iter{n}}(\hat{\bm{m}}) = \left\{ \bm{x} \notin \mathcal{H} \,:\, \lambda^{\iter{n}}(\bm{x}) > \lambda^{\iter{n}}(\hat{\bm{m}})\right\}.
\]
We then sort all candidate midpoints in the set $\bm{X}^{\iter{n}}_\psi$, by decreasing $\lambda^{\iter{n}}$ values:

{\small
\[
    \bm{M}^{\iter{n}} \triangleq\left\{\bm{m}^{\iter{n}}_1,\ldots,\bm{m}^{\iter{n}}_{M}\right\}:\lambda^{\iter{n}}(\bm{m}^{\iter{n}}_k) \geq \lambda^{\iter{n}}(\bm{m}^{\iter{n}}_{k+1})
\]
}

\noindent where $\bm{m}^{\iter{n}}_k=\bm{x}^{(i,j)}_\psi$ for some $\bm{x}^{(i)},\bm{x}^{(j)} \in \bm{X}^{\iter{n}}$ (see \eqref{eqn:explore-set}) and $M=\dfrac{n(n-1)}{2}$ is the total number of segments at iteration $n$. In Mode~$\psi$, the first-ranking candidate $\bm{m}^{\iter{n}}_1$ is taken for sampling. Its uncertainty is necessarily larger than $\lambda^{\iter{n}}(\hat{\bm{m}})$ (otherwise $\hat{\bm{m}}$ would have been ranked higher), i.e. $\bm{m}^{\iter{n}}_1\in\mathcal{R}^{\iter{n}}(\hat{\bm{m}})$. Assume, for the purpose of contradiction, that points in $\mathcal{H}$ are never sampled. 

\noindent At iteration $n+1$, the midpoint $\bm{m}^{\iter{n}}_1\notin \mathcal{H}$ is thus sampled and its uncertainty $\lambda^{\iter{n+1}}(\bm{m}^{\iter{n}}_1)$ becomes zero. Since this point does not belong to $\mathcal{H}$, the value of $h$ in \eqref{eqn:midpoint-h} is still valid, while the uncertainty bound \eqref{eqn:midpoint-proof} may either be the same, or increase in case $\gamma^{\iter{n+1}}>\gamma^{\iter{n}}$ (see \eqref{eqn:gamma-update}).  Moreover, for any point $\bm{x}\in\mathcal{X}$ we have
\begin{equation}\label{eqn:lemma4-1}
\lambda^{\iter{n+1}}(\bm{x})\leq2\mu\gamma^{\iter{n+1}}\|\bm{x}-\bm{m}^{\iter{n}}_1\|
\end{equation}
as shown by direct application of \eqref{eqn:lower-bounds}-\eqref{eqn:sm-uncertainty}. Consider now the ball $\mathcal{B}(\bm{m}^{\iter{n}}_1,r_\lambda)$, where
\begin{equation}\label{eqn:lemma4-2}
    r_\lambda = \left(1 - \frac{1}{\mu}\right) h.
\end{equation}
On the basis of \eqref{eqn:midpoint-proof}, \eqref{eqn:lemma4-1} and \eqref{eqn:lemma4-2} we have that:

\begin{multline}
\label{eqn:ineq-ball}
    \forall \bm{x} \in \mathcal{B}(\bm{m}^{\iter{n}}_1,r_\lambda),\; \lambda^{\iter{n+1}}(\bm{x}) \leq  2\mu\gamma^{\iter{n+1}}r_\lambda = \\ 2(\mu-1)\gamma^{\iter{n+1}} h\leq \lambda^{\iter{n+1}}(\hat{\bm{m}}).
\end{multline}

\noindent Hence,
\[
\mathcal{R}^{\iter{n+1}}(\hat{\bm{m}}) = \mathcal{R}^{\iter{n}}(\hat{\bm{m}}) \setminus \left( \mathcal{R}^{\iter{n}}(\hat{\bm{m}}) \cap \mathcal{B}(\bm{m}^{\iter{n}}_1, r_\lambda) \right).
\]

\noindent This means that the volume of set $\mathcal{R}^{\iter{n}}(\hat{\bm{m}})$ diminishes by a finite quantity at each iteration, since the value of $h$ does not depend on $n$; and becomes null after finite iterations, since Mode $\psi$ is assumed to be undertaken infinitely often. Denote with $\bar{n}-1<\infty$ the iteration at which this happens. Then, no midpoint outside the set $\mathcal{H}$ would have uncertainty larger than $\lambda^{\iter{\bar{n}}}(\hat{\bm{m}})$, so that a midpoint 

\[
    \left\{ \idxpt{x}{i,j}{\psi} \in \bm{X}^{\iter{\bar{n}}}_\psi \cap \mathcal{H} ~:~ \lambda^{\iter{\bar{n}}}(\idxpt{x}{i,j}{\psi}) \geq \lambda^{\iter{\bar{n}}}(\hat{\bm{m}})\right\}
\]

\noindent will be ranked first in $\bm{M}^{\iter{\bar{n}-1}}$ and sampled at iteration $\bar{n}$, thus falling in contradiction.\\
We thus demonstrated that, in finite iterations, a point inside $\mathcal{H}$ is sampled. Moreover, denoting such a point as $\idxpt{x}{\bar{n}}{}$, by construction we have 
\begin{equation}
\label{eqn:radius-smaller}
\| \idxpt{x}{\bar{n}}{} - \bar{\bm{x}}^{\iter{n}} \| < \| \tilde{\bm{x}}^{\iter{n}} - \bar{\bm{x}}^{\iter{n}} \|.
\end{equation}
Then, two cases may occur: either $\idxpt{x}{\bar{n}}{}$ is closer to $\hat{\bm{x}}$ than $\bar{\bm{x}}^{\iter{n}}$, or not. In the former case, repeating the process from \eqref{eqn:lemma-nearest-sample} we'll have that $\bar{\bm{x}}^{\iter{\bar{n}}}=\idxpt{x}{\bar{n}}{}$ and $\tilde{\bm{x}}^{\iter{\bar{n}}}=\bar{\bm{x}}^{\iter{n}}$. In the latter case, $\bar{\bm{x}}^{\iter{\bar{n}}}=\bar{\bm{x}}^{\iter{n}}$ and $\tilde{\bm{x}}^{\iter{\bar{n}}}=\idxpt{x}{\bar{n}}{}$. In both cases, the radius of the half-hyper-ball $\mathcal{H}$ shrinks because of \eqref{eqn:radius-smaller}. By induction we thus have:

\begin{equation}
\label{eqn:half-ball-shrinks}
    \lim\limits_{n\rightarrow\infty}\| \tilde{\bm{x}}^{\iter{n}} - \bar{\bm{x}}^{\iter{n}} \| =0.
\end{equation}

Next, still considering $\bar{\bm{x}}^{\iter{n}}$ as the nearest sample to $\hat{\bm{x}}$ \eqref{eqn:lemma-nearest-sample}, we are going to show that there exists a radius $a>0$ such that 
\begin{equation}
\label{eqn:finite-radius-nearer}
\| \bm{x} - \hat{\bm{x}} \| ~<~ \| \bar{\bm{x}}^{\iter{n}} - \hat{\bm{x}} \|, \forall \bm{x}\in\mathcal{Q}\cap\mathcal{B}(\bar{\bm{x}}^{\iter{n}},a)
\end{equation}

To this end, consider any point  $\bm{x}\in\mathcal{Q}$ and assume, for the sake of contradiction, that \eqref{eqn:finite-radius-nearer} never holds. This would mean that the nearest point to $\hat{\bm{x}}$ on the segment connecting $\bar{\bm{x}}^{\iter{n}}$ to $\bm{x}$ is always $\bar{\bm{x}}^{\iter{n}}$, no matter how close the two points are. In turn, this would imply that:
\begin{equation}
\label{eqn:directional-ball}
    \left.\frac{d}{da} \|  ((1-a)\bar{\bm{x}}^{\iter{n}} + a\bm{x}) - \hat{\bm{x}} \| \right|_{a=0}\geq 0, 
\end{equation}

\noindent i.e.,
\begin{align*}
    \frac{1}{\| \bar{\bm{x}}^{\iter{n}} - \hat{\bm{x}} \|} (\bar{\bm{x}}^{\iter{n}} - \hat{\bm{x}})^\top (\bm{x}-\bar{\bm{x}}^{\iter{n}}) &\geq 0 \\
    (\hat{\bm{x}} - \bar{\bm{x}}^{\iter{n}})^\top (\bm{x}-\bar{\bm{x}}^{\iter{n}}) &\leq 0
\end{align*}

\noindent which would imply that $\bm{x}\notin\mathcal{Q}$, thus falling in a contradiction. Hence, the directional derivative \eqref{eqn:directional-ball} must be negative for all $\bm{x}\in\mathcal{Q}$, meaning that there exists a scalar $a>0$ such that in the neighborhood $\mathcal{Q}\cap\mathcal{B}(\bar{\bm{x}}^{\iter{n}},a)$ the property \eqref{eqn:finite-radius-nearer} holds.

Combining \eqref{eqn:half-ball-shrinks} and \eqref{eqn:finite-radius-nearer}, we obtain that, after a finite number $k$ of iterations, the following condition holds:
\[
    \| \bar{\bm{x}}^{\iter{n+k}} - \hat{\bm{x}} \| ~<~ \| \bar{\bm{x}}^{\iter{n}} - \hat{\bm{x}} \|.
\]
\noindent This implies that
\[ 
    \lim\limits_{n \rightarrow \infty}\| \bar{\bm{x}}^{\iter{n}} - \hat{\bm{x}} \| = 0,
\]
\noindent and, for any $\sigma>0$, we finally have
\[
    \exists n_\sigma: \min_{\bm{x}^{(i)} \in \bm{X}^{\iter{n_\sigma}}} \| \bm{x}^{(i)} - \hat{\bm{x}} \| < \sigma
\]
\noindent thus proving the lemma.
\end{proof}

\begin{theorem}
Let Assumptions \ref{ass:lipschitz} and \ref{assumption:assume-corners-sampled} hold. Then,
$\forall \epsilon>0,\;\exists n_\epsilon<\infty: z^{*\iter{n_\epsilon}} \leq z^* + \epsilon$.
\end{theorem}

\begin{proof}
Consider a point $\hat{\bm{x}}\in\mathcal{X}^*$ (see \eqref{eqn:global-minimzer}) and take $\sigma = \dfrac{\epsilon}{\gamma_o}$. For any $n$, denote
\[
\bar{\bm{x}}^{\iter{n}}=\arg\min\limits_{\bm{x}^{(i)} \in \bm{X}^{\iter{n}}} \| \hat{\bm{x}} - \bm{x}^{(i)} \|
\]
In virtue of Lemma \ref{lemma:exploit-will-fail}, the exploration mode will be called infinitely often. 
Now, by applying Lemma \ref{lemma:dense-pts}, we have that $\exists n_\sigma < \infty :  \| \hat{\bm{x}} -\bar{\bm{x}}^{\iter{n_\sigma}}\| ~<~ \sigma=\dfrac{\epsilon}{\gamma_o}$. 
Then, in virtue of Assumption \ref{ass:lipschitz} we have:

\begin{multline}
    f_o(\bm{x}^{*\iter{n_\sigma}})-f_o(\hat{\bm{x}})=z^{*\iter{n}}-z^*\leq \\ f_o(\bar{\bm{x}}^{\iter{n_\sigma}})-f_o(\hat{\bm{x}})\leq  \gamma_o\|\hat{\bm{x}} - \bar{\bm{x}}^{\iter{n_\sigma}} \| ~< \epsilon
\end{multline}

Thus proving the result with $n_\epsilon=n_\sigma$.
\end{proof}

\subsection{Optimality gap}\label{subsec:optimality-gap}
Considering a run of the SMGO algorithm with $N$ iterations, we now analyze the gap between the best sampled value and the actual global minimum,

\[
    \delta^{\iter{N}} = z^{\iter{N}*} - f_o(x^*)
\]

\noindent denoted as the \textit{optimality gap}. In calculating such a difference, the SM-based guarantees $\underline{z}^{\iter{N}}(\bm{x}) \leq f_o(\bm{x}) \leq \overline{z}^{\iter{N}}(\bm{x})$ can be utilized in the case of a known Lipschitz constant $\gamma_o$. Hence, the upper bound of the optimality gap $\overline{\delta}^{\iter{N}}$ can be considered as the difference between $z^{\iter{N}*}$ and the minimum lower bound $\underline{z}^{\iter{N}}$ defined in \eqref{eqn:lower-bounds} generated by data set $\bm{X}^{\iter{N}}$, i.e.

\begin{align}
  \overline{\delta}^{\iter{N}} &= \max_{\bm{x} \in \mathcal{X}} \left( z^{\iter{N}*} - \underline{z}^{\iter{N}}(\bm{x}) \right) \nonumber\\
  ~&= z^{\iter{N}*} - \min_{\bm{x} \in \mathcal{X}} \underline{z}^{\iter{N}}(\bm{x})
  \label{eqn:subopt-gap}
\end{align}

The calculation of $\min\limits_{\bm{x} \in \mathcal{X}} \underline{z}^{\iter{N}}(\bm{x})$ in \eqref{eqn:subopt-gap} entails calculation of the hyperbolic Voronoi cells and of the lower bounds at their vertices. Denoting with $\bm{V}^{\iter{N}}$ the set of such vertices generated from $\bm{X}^{\iter{N}}$, \eqref{eqn:subopt-gap} can be expressed as

\begin{equation}
  \label{eqn:subopt-gap-vertex}
  \overline{\delta}^{\iter{N}} = z^{\iter{N}*} - \min_{\bm{v} \in \bm{V}^{\iter{N}}}\underline{z}^{\iter{N}}(\bm{v}).
\end{equation}

The acquisition of HVD vertices using the samples can be approached as described in~\cite{Milanese2007}, which is referred for more information. The calculation of HVD cells and vertices entail computations with exponential complexity w.r.t. $D$~\cite{Milanese2007}. However, an interested user can obtain an upper bound on the optimality gap by carrying out \eqref{eqn:subopt-gap-vertex}.

\section{Computational Aspects}
\label{sec:implement-notes}
This section presents an analysis of the computational aspects of the proposed algorithm and possible improvements. The first aspect discussed is on ensuring the coverage of $\mathcal{X}$, followed by a discussion on the computational complexity.

\subsection{On the coverage of the search set}\label{subsec:coverage}

The selection method for candidate points described in Section \ref{sec:smgo-algo}  allows SMGO to only explore the volume inside the convex hull of the sampled points in $\bm{X}^{\iter{n_0}}$. Hence, in order to ensure exploration throughout all the search set $\mathcal{X}$, all its vertices $\idxpt{c}{v}{}$ should also be considered in the candidate points generation, as stated in Assumption \ref{assumption:assume-corners-sampled}. One approach is to first evaluate $f_o$ on all $\idxpt{c}{v}{}$ before actually starting with exploitation/exploration. However, in the case of a hyperrectangle, this entails an exponential number of long function calls w.r.t. dimensionality $D$.

\begin{figure}[!t]
	\centering
	\includegraphics[width=\columnwidth]{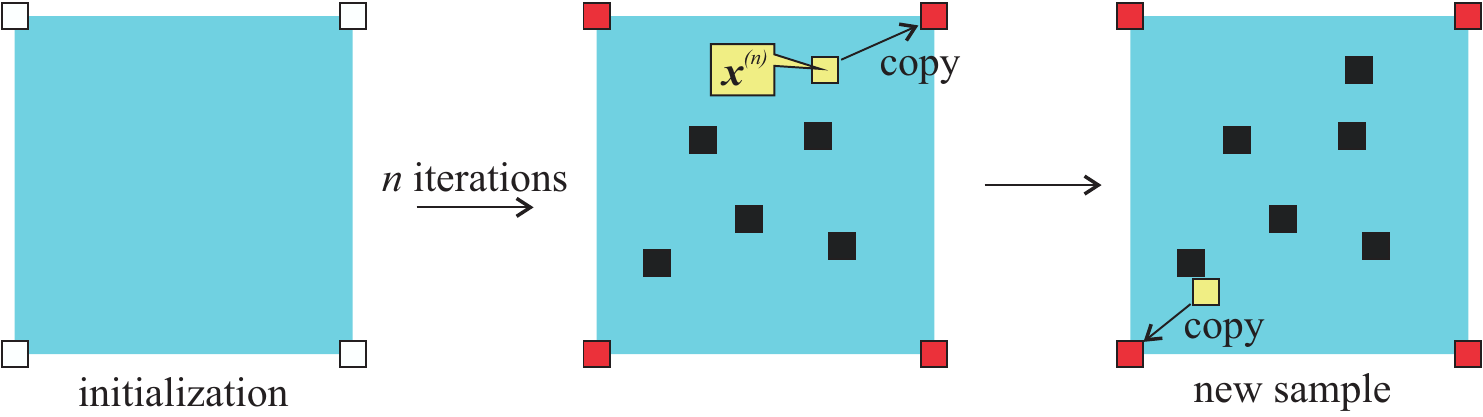}
	\caption{Corner mirroring of nearest sample point}
	\label{fig:coverage}
\end{figure}

To address this issue, a simple approach is proposed in order to ensure coverage in $\mathcal{X}$ without explicitly calling $f_o$ to evaluate the corners. Moreover, with this approach Assumption \ref{assumption:assume-corners-sampled} is not needed anymore for the theoretical results to hold. At each iteration $n$, the midpoints between any sample and all the vertices of $\mathcal{X}$ are included as candidate exploration samples, even if the vertices do not belong to $\bm{X}^{\iter{n}}$. Moreover, for the purpose of computing the uncertainty of each midpoint, the cost function value $f_o(\idxpt{c}{v}{})$ at the $v$-th vertex $\idxpt{c}{v}{}$ is estimated as: 

\begin{equation}\label{eqn:vertex-trick}
  f_o(\idxpt{c}{v}{}) \approx z^{(w)},\;\text{where\;}
w = \arg\min\limits_{k=1, \ldots, n} \| \bm{x}^{(k)} - \idxpt{c}{v}{} \|
\end{equation}
i.e., equal to the cost function value at the nearest sampled point.
The concept of this strategy is shown in Fig.~\ref{fig:coverage}. Thanks to Assumption \ref{ass:lipschitz} and Lemma \ref{lemma:dense-pts}, with an increasing number of sample points in $\mathcal{X}$  the estimate \eqref{eqn:vertex-trick} approximates more and more accurately the true function value at the vertex.

\subsection{On the computational burden}
\label{subsec:segment-gen}

There are three main routines within an iteration of SMGO: Lipschitz constant update, exploitation by Mode~$\theta$ routine, and eventually, if the expected improvement threshold is not achieved, exploration by Mode~$\psi$. The worst-case SMGO complexity analysis assumes that each iteration eventually proceeds to Mode~$\psi$.
 
The first routine compares the incoming $n$-th sample with $n-1$ existing ones, resulting in $\mathcal{O}(n)$ complexity.
 
The Modes~$\theta$ and $\psi$ can be designed as iterative routines as well to significantly improve the computational burden, as discussed in the following.

The exploitation Mode~$\theta$ comprises the selection of $n-1$ candidate points $\bm{x}^{(i)}_\theta$ for each segment $l_{n^*,i}$ from $\bm{x}^{*\iter{n}}$ to all the other sampled points according to \eqref{eqn:exploit1-const}-\eqref{eqn:exploit1-loc}, and the calculation of $\underline{z}^{\iter{n}}$ for each $\bm{x}^{(i)}_\theta$ using \eqref{eqn:lower-bounds}. Both are $\mathcal{O}(n)$ operations, leading to a compounded complexity of $\mathcal{O}(n^2)$. However, when $\gamma^{\iter{n}}$ and $\bm{x}^{*\iter{n}}$ do not change from iteration $n-1$ to $n$, a cache can store $\bm{x}^{(i)}_\theta$ locations and corresponding $\underline{z}^{\iter{n}}$. From \eqref{eqn:lower-bounds}, it is convenient to save $z^{(\tilde{k})}$ and $\tilde{\Delta}^{(i)} \triangleq \mu\gamma^{\iter{n}}\| \bm{x}^{(i)}_\theta - \bm{x}^{(\tilde{k})}\|$, where 

\[ \tilde{k} = \argmax_{k \in [1 \ldots n]} \left( z^{(k)}-\mu\gamma^{\iter{n}}\| \bm{x}^{(i)}_\theta - \bm{x}^{(k)}\| \right).\] 

\noindent In practice, information are saved\footnote{It can be understood that $z^{(\tilde{k})}$ is the value at the tip, while $\tilde{\Delta}^{(i)}$ is the height of the hypercone.} about the hypercone that actually generated $\underline{z}^{\iter{n}}$ for each $\bm{x}^{(i)}_\theta$. As long as $\bm{x}^{*\iter{n}}$ stays the same, previously existing $\bm{x}^{(i)}_\theta$ locations do not need to be computed again.

Furthermore, calculating $\underline{z}^{\iter{n}}$ for an existing exploitation point $\bm{x}^{(i)}_\theta$ reduces to

\begin{equation}
\label{eqn:theta-iter-lower-bound}
\underline{z}^{\iter{n}}(\bm{x}^{(i)}_\theta) = \max\left(\underline{z}^{\iter{n-1}}(\bm{x}^{(i)}_\theta), \underline{z}^{(i)}_{new} \right)
\end{equation}

\noindent where

\[
    \underline{z}^{(i)}_{new} = z^{(n)} - \mu\gamma^{\iter{n}} \| \bm{x}^{(i)}_\theta - \bm{x}^{(n)} \|
\]

\noindent for an incoming sample $(\bm{x}^{(n)}, z^{(n)})$. 

Note that when $\underline{z}^{\iter{n}}(\bm{x}^{(i)}_\theta)$ changes due to \eqref{eqn:theta-iter-lower-bound}, $z^{(\tilde{k})} = z^{(n)}$ and $\tilde{\Delta}^{(i)} = \mu\gamma^{\iter{n}}\| \bm{x}^{(i)}_\theta - \bm{x}^{(n)}\|$. As a result, computing $\underline{z}^{\iter{n}}$ for all existing candidate points is $\mathcal{O}(n)$, i.e. $\mathcal{O}(1)$ for each. Now, consider the situation when $\gamma^{\iter{n}}$ is updated, assuming that $\bm{x}^{*\iter{n}}$ and $\tilde{k}$ do not change. This implies rescaling $\tilde{\Delta}^{(i)}$ such that

\begin{equation}
    \label{eqn:theta-iter-gamma-changed}
    \tilde{\Delta}^{(i)}_{\iter{n}} = \frac{\gamma^{\iter{n}}}{\gamma^{\iter{n-1}}} \tilde{\Delta}^{(i)}_{\iter{n-1}}
\end{equation}

\noindent and updating $\underline{z}^{\iter{n}}$ with the new sample $(\bm{x}^{(n)}, z^{(n)})$ using \eqref{eqn:theta-iter-lower-bound}, still resulting in $\mathcal{O}(1)$ per candidate point. This implies a Mode~$\theta$ aggregate complexity of $\mathcal{O}(n)$.

Finally, one new candidate point $\bm{x}^{(n)}_\theta$ is evaluated using \eqref{eqn:lower-bounds}, which is $\mathcal{O}(n)$. Hence, as long as $\bm{x}^{*\iter{n}}$ does not change, the introduction of a cache reduces Mode~$\theta$ to $\mathcal{O}(n)$. However, when $\bm{x}^{*\iter{n}}$ changes, i.e. $\bm{x}^{*\iter{n}} = \bm{x}^{(n)}$, computing all $\bm{x}^{(i)}_\theta$ locations must be repeated from the start, leading to $\mathcal{O}(n^2)$ worst-case complexity for Mode~$\theta$ as stated before.

As defined, Mode~$\psi$ generates $\dfrac{n(n-1)}{2}$ candidate points (midpoints among all existing samples) at each iteration. Furthermore, the calculation of $\lambda^{\iter{n}}$ for each candidate point is $\mathcal{O}(n)$, resulting in $\mathcal{O}(n^3)$ for each exploration iteration.

Note that behavior in Mode $\psi$ does not depend on the current best point, because $n$ new segments are added to the search set at iteration $n$. A cache solution can be introduced also in this mode, updating $\lambda^{\iter{n}}$-related values iteratively. 

For each existing $\bm{x}^{(i,j)}_\psi$, the values to be saved are $z^{(\hat{k})}$ and $\hat{\Delta}^{(i,j)}$ (for $\overline{z}^{\iter{n}}$), and $z^{(\check{k})}$ and $\check{\Delta}^{(i,j)}$ (for $\underline{z}^{\iter{n}}$), where $\hat{k}$ and $\check{k}$ are defined as

\[\hat{k} = \argmin_{k \in [1 \ldots n]} \left( z^{(k)}+\mu\gamma^{\iter{n}}\| \bm{x}^{(i,j)}_\psi - \bm{x}^{(k)}\| \right),\] 

\[\check{k} = \argmax_{k \in [1 \ldots n]} \left( z^{(k)}-\mu\gamma^{\iter{n}}\| \bm{x}^{(i,j)}_\psi - \bm{x}^{(k)}\| \right).\] 

The iterative updates are performed using eq. \eqref{eqn:theta-iter-lower-bound} and \eqref{eqn:theta-iter-gamma-changed}. This then results in $\mathcal{O}(1)$ for each existing midpoint, and $\mathcal{O}(n)$ for each new midpoint generated by new segments from incoming $\bm{x}^{(n)}$ to the $n-1$ existing points. This results in $\mathcal{O}(n^2)$ for Mode~$\psi$, and in turn, a $\mathcal{O}(n^2)$ worst-case complexity for SMGO.

\begin{remark}

Note that when $\gamma^{\iter{n}}$ is updated (which always implies a growth of its value), the proposed iterative implementation might result in more conservative bounds $\overline{z}^{\iter{n}}$ and $\underline{z}^{\iter{n}}$ than if repeatedly recalculated at every iteration, i.e. lower $\underline{z}^{\iter{n}}$ for Mode~$\theta$, and larger $\lambda^{\iter{n}}$ in Mode~$\psi$. This occurs because with higher $\gamma^{\iter{n}}$-values the HVD layout increasingly approaches that of a non-hyperbolic Voronoi tessellation, i.e. the cones providing the tightest bounds at any $\idxpt{x}{i,j}{\psi}$ could change when $\gamma^{\iter{n+1}} \gg \gamma^{\iter{n}}$, and $\itrval{z}{\hat{k}}{}$ and $\itrval{z}{\check{k}}{}$ would tend to the nearest sample w.r.t. said $\idxpt{x}{i,j}{\psi}$. These effects are not accounted for by the iterative update procedure, which assumes that the bound-generating cones remain the same. However, the theoretical properties of SMGO are still valid with the iterative implementation. In fact, all the arguments for Lemma~\ref{lemma:exploit-will-fail} are valid also with iteratively-computed lower bounds. Regarding the exploration mode, inequality (\ref{eqn:min-lambda}) in Lemma~\ref{lemma:theoretical-min-lambda} is consistent with the more conservative iteratively-computed $\gamma^{\iter{n}}$. Furthermore, when sampling a new point $\itrpt{m}{n}{1}$, for which $\lambda^{\iter{n+1}}(\itrpt{m}{n}{1})$ is exactly computed, (\ref{eqn:ineq-ball}) still holds. Hence, all arguments laid out in Lemma~\ref{lemma:dense-pts} still hold.

\end{remark}

\section{Performance Test Results}
\label{sec:perf-test}
In this section, the  proposed SMGO algorithm is evaluated on well-known benchmarks, usually employed in black-box optimization, comparing its performance with representative Lipschitz-based methods: DIRECT and AdaLIPO. In addition, the Bayesian optimization approach (BayesianO) is also considered in the tests. The results are discussed in terms of iteration-based optimization performance and computational times. All scripts used for this section are available at \verb|https://github.com/lorenzosabugjr/smgo|.

\subsection{Test parameters}
Seven test functions of varying structure and dimensions are chosen for this test, covering a variety of characteristics relevant to black-box optimization benchmarks. Their search bounds, optimal values and relevant features are summarized in Table~\ref{table:test-fns}. A reasonable assumption in real-life applications is that the black-box function has no available optimal value. Hence, each algorithm is given a budget of $N=500$ long function evaluations for each optimization run. The optimal results after the evaluation budget are then measured and compared among the different algorithms, in what is referred as a fixed-cost/budget comparison~\cite{Beiranvand2017}. All computations were performed in MATLAB 2020b on a system with AMD Ryzen 9 3900X (3.80 GHz) and 32 GB RAM.

To assess the sensitivity of each algorithm to different initial information and randomized decisions, 100 independent runs (trials) were performed. For each trial, the same randomly-generated starting point is given to AdaLIPO, SMGO, and BayesianO. On the other hand, DIRECT considers a fixed set of search points, not allowing to randomize the starting samples. 

SMGO is implemented as described in Algorithm \ref{algo:smgo} with parameters $\alpha=0.001$ and $\mu=1.025$. AdaLIPO is implemented as described in \cite{Malherbe2017} with  parameter $p=0.1$. DIRECT is implemented as described in \cite{Finkel2006}, not requiring any tuning parameter. Finally, the algorithm available in the Matlab Statistics and Machine Learning Toolbox is employed for the BayesianO implementation, using the provided default parameters.

\begin{table*}[!t]
    \footnotesize
    \begin{center}
    	\begin{tabular}{|c|p{2.0in}|c|c|p{0.75in}|}
    		\hline
    		\textbf{Description} & \textbf{Function definition} $f(\bm{x})$ & \textbf{Bounds} & \textbf{$z^*$} & \textbf{Features}\\ \hline
    		Rosenbrock & $\sum_{i=1}^{D}[100 (x_{i+1} - x_i^2)^ 2 + (1 - x_i)^2]$ & $-40 \leq x_i \leq 5$ & 0.0 & unimodal, non-separable \\ \hline
    		Styblinski-Tang & $\frac{1}{2} \sum_{i=1}^D \left( x_i^4 - 16x_i^2 + 5x_i \right)$ & $-5 \leq x_i \leq 5$ & $-39.166 D$ & multimodal, separable \\ \hline
    		Deb's \#1 & $\frac{1}{D} \sum_{i=1}^D\textrm{sin}^6(5\pi x_i)$ & $-1 \leq x_i \leq 1$ & $-1.0$ & multimodal, separable, numerous global minima \\ \hline
    		Deb's \#2 & $\frac{1}{D} \sum_{i=1}^D\textrm{sin}^6[5\pi (x_i^{3/4} -0.05)]$ & $0 \leq x_i \leq 150$ & $-1.0$ & multimodal, separable, numerous global minima \\ \hline
    		Schwefel & $- \sum_{i=1}^D x_i \textrm{sin}\left(\sqrt{|x_i|}\right)$ & $-500 \leq x_i \leq 500$ & $-418.982 D$ & multimodal, separable, numerous local minima \\ \hline
    		Salomon & $1-\textrm{cos}\left(2\pi\sqrt{\sum_{i=1}^D x_i^2}\right) + 0.1\sqrt{\sum_{i=1}^D x_i^2}$ & $-40 \leq x_i \leq 70$ & $0$ & multimodal, non-separable \\ \hline
    		Brown & $\sum_{i=1}^{D-1} (x_i^2)^{(x_{i+1}^2+1)} + (x_{i+1}^2)^{(x_i^2+1)}$ & $-1 \leq x_i \leq 4$ & $0$ & multimodal, non-separable \\ \hline
    	\end{tabular}
        \caption{Test functions used for comparative tests}
    	\label{table:test-fns}
    \end{center}
\end{table*}

\begin{figure*}[!t]
	\centering
 	\includegraphics[width=0.625\textwidth]{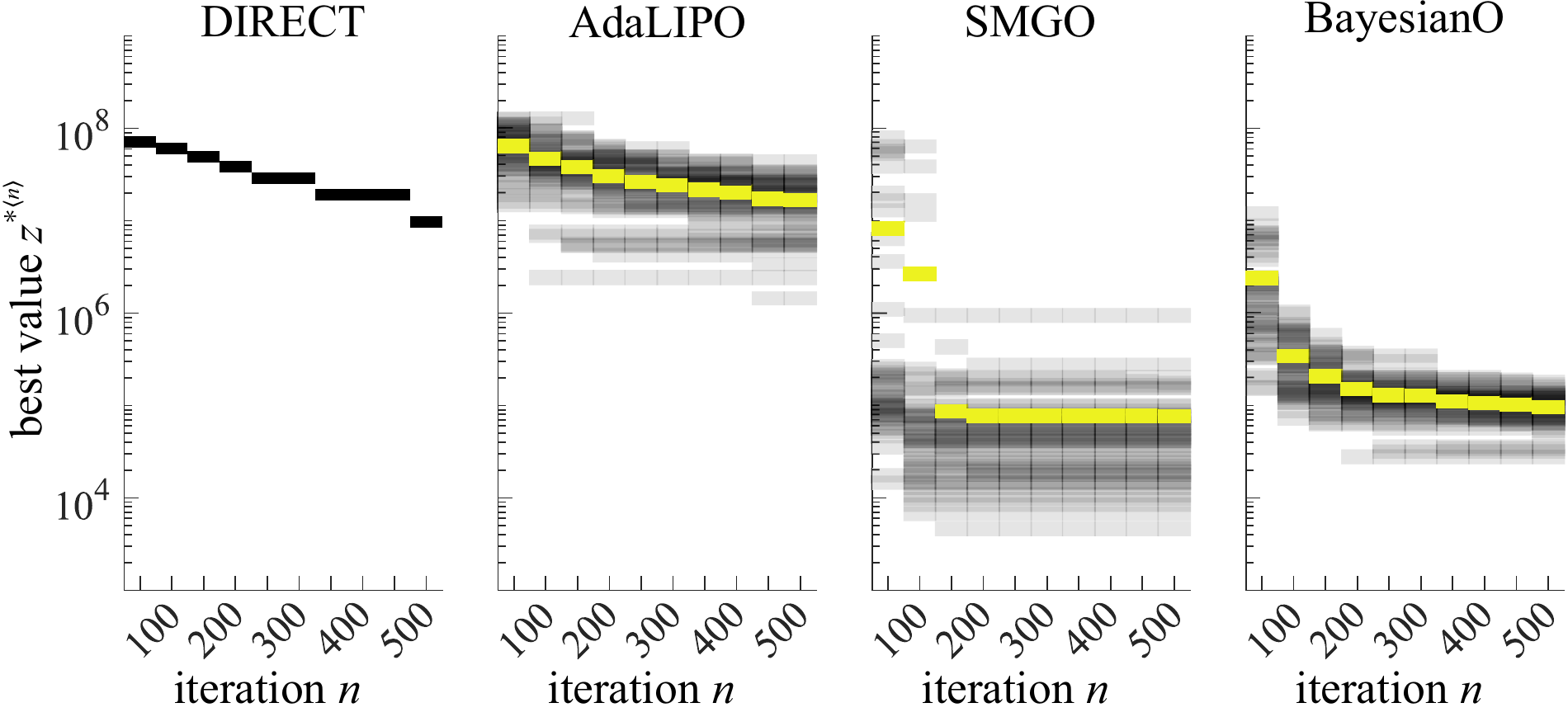}
	\caption{Optimal value distribution w.r.t. iterations, \\10D Rosenbrock function (log scale)}
	\label{fig:results-rosen-10d}
\end{figure*}

\begin{figure*}[!t]
	\centering
	\includegraphics[width=0.625\textwidth]{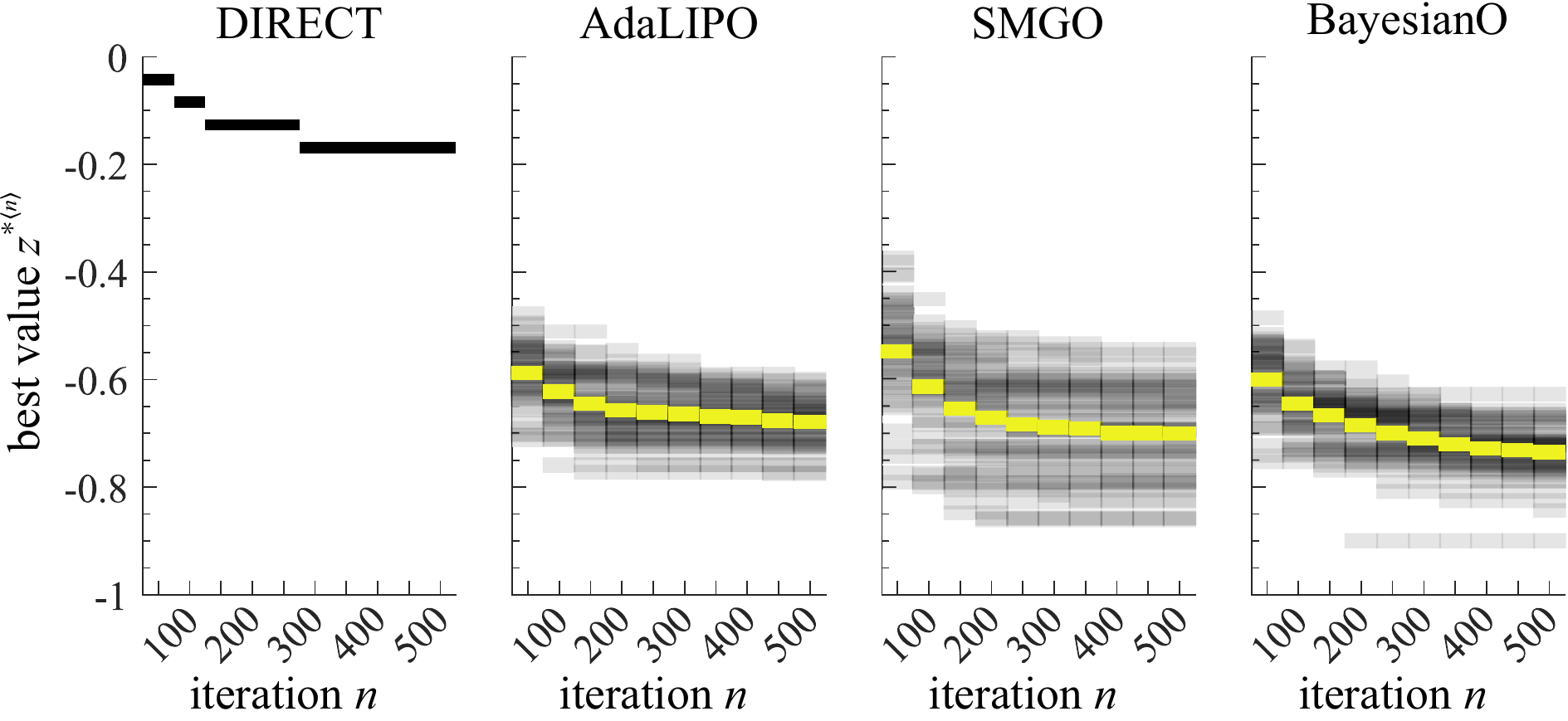}
	\caption{Optimal value distribution w.r.t. iterations, \\10D Deb's function \#1}
	\label{fig:results-deb1-10d}
\end{figure*}

\subsection{Iteration-based performance}

The evolution of the best sample $z^{*\iter{n}}$  along the iterations ($n\in [50, 500]$) for each algorithm is shown in Figs.~\ref{fig:results-rosen-10d} and \ref{fig:results-deb1-10d} for the Rosenbrock and Deb's functions, respectively. The graphs show the distribution of the best sampled cost for the 100 independent runs. Note that the results of the DIRECT algorithm do not show dispersion due to the fact that it is a deterministic search method.  DIRECT follows a batch sampling strategy, hence the graphs show for each $n$ value, the best results for the previously-sampled batches.

The Rosenbrock function is a standard test case for gradient-based algorithms, with a unique global minimum (unimodal) for $D=2$; however it has 2 minima for $D=4\sim30$~\cite{Shang2006}. The location of $\bm{x}^*$ is at $(1, 1,\ldots, 1)$, hence the classical search bounds $-10 \leq x_i \leq 10$ would be advantageous for DIRECT, which invariably samples the center of the search hyper-box. For this reason, the search bounds are changed to move $\bm{x}^*$ away from the center. 

As seen in Fig.~\ref{fig:results-rosen-10d}, DIRECT has only improved marginally the initial sampled cost after the 500 iterations. Also AdaLIPO achieved  marginal improvements after 500 iterations. This might be associated to the completely randomized exploration and exploitation decisions. On the other hand, both SMGO and BayesianO achieved large improvements during the first 100 iterations, while  slight enhancements are observed during the final 400 iterations. SMGO exhibits a wider distribution among different runs and, on average, achieves a better result after 500 iterations.

In the case of Deb's function \#1, as shown in Fig.~\ref{fig:results-deb1-10d}, the presence of numerous global minima  causes a bad performance for DIRECT, which barely improved its initial best cost. On the other hand, AdaLIPO, BayesianO, and SMGO have resulted in similar results, achieving very similar average costs after the iteration limit, while BayesianO shows more concentrated values among runs. In this case the randomized mechanism of AdaLIPO allows it to discover low-cost regions to improve the best sample. Furthermore, Bayesian's usage of a response surface, and SMGO's geometry-based candidate points derivation have become useful in looking for low-cost regions in the search space.

A summary (mean best values) of the algorithms' outcomes after 500 iterations are shown in Table~\ref{table:performance-results}. Note that the results tabulated for DIRECT were achieved at iterations which are a bit larger than the pre-defined maximum $N=500$, owing to its batch-based samples generation. In most test cases, BayesianO achieved the best average results, while SMGO and DIRECT attained the best average performance in 2 cases each. On the other hand, AdaLIPO achieved fair results, even if it ranked second-place or worse for most of the tested functions.

Furthermore, for each test function (with the same starting points across all algorithms), the number of trials (out of 100) resulting in best results for each respective algorithm is shown in Table~\ref{table:best-results}. As with the average performance, BayesianO achieved the most best trials in 9 out of the 13 tested functions. Furthermore, in Deb's \#1 (5D), Schwefel (5D and 10D), and Salomon (5D and 10D), it resulted in best results in almost all of the given trials. On the other hand, SMGO had the most best results in Rosenbrock and Deb's \#2 (10D), and had the second largest count of best trials in Deb's \#1 (10D) and Deb's \#2 (5D). AdaLIPO registered the second count of best trials in Deb's \#2 (10D), and had no best trials for 4 functions. Lastly, DIRECT had no best trials for 8 out of 13 test functions, but took almost all the best results in the Brown function (5D and 10D). 

We have also compared the outcomes of SMGO with the competitor methods using the Wilcoxon and Kruskal-Wallis non-parametric statistical tests, pairing the results by the shared random initial point. Both test were performed with 5\% significance level and with the null hypothesis that SMGO results are statistically similar to those from the respective competitor method. The two tests agreed that SMGO results were statistically similar to AdaLIPO for Styblinski-Tang (5D), Deb's \#2 (5D, 10D), Schwefel (10D), and Brown (10D). Furthermore, similarities were found between SMGO and BayesianO for Deb's \#2 (10D). Given these results, SMGO is found to have competitive results compared with the considered global optimization algorithms on the test functions employed in the analysis.

\begin{table*}[!t]
\footnotesize
\begin{center}
	\begin{tabular}{|c|c|r|r|r|r|}
		\hline
		\textbf{Test function} & $D$ & \textbf{DIRECT} & \textbf{AdaLIPO} & \textbf{SMGO} & \textbf{BayesianO} \\ \hline
		                  Rosenbrock     & 10 &  1.17\,E\,+5 &  1.77\,E\,+7  &  \textbf{8.63\,E\,+4} &  9.35\,E\,+4          \\ \hline
		\multirow{2}{*}{Styblinski-Tang} & 5  &  -1.95\,E\,+2   &  -1.59\,E\,+2  &  -1.58\,E\,+2           &  \textbf{-1.95\,E\,+2}  \\ \cline{2-6}
		                                 & 10 &  -3.19\,E\,+2   &  -2.61\,E\,+2  &  -2.96\,E\,+2           &  \textbf{-3.33\,E\,+2}  \\ \hline
		\multirow{2}{*}{Deb's \#1}       & 5  &  -5.29\,E\,$-$1     &  -8.29\,E\,$-$1    &  -8.07\,E\,$-$1             &  \textbf{-9.68\,E\,$-$1}    \\ \cline{2-6}
		                                 & 10 &  -2.11\,E\,$-$1     &  -6.74\,E\,$-$1    &  -6.97\,E\,$-$1             &  \textbf{-7.35\,E\,$-$1}    \\ \hline
		\multirow{2}{*}{Deb's \#2}       & 5  &  -5.97\,E\,$-$1     &  -8.32\,E\,$-$1    &  -8.33\,E\,$-$1             &  \textbf{-8.59\,E\,$-$1}    \\ \cline{2-6}
		                                 & 10 &  -3.98\,E\,$-$1     &  -6.72\,E\,$-$1    &  \textbf{-6.81\,E\,$-$1}    &  -6.77\,E\,$-$1             \\ \hline
		\multirow{2}{*}{Schwefel}        & 5  &  -1.47\,E\,+3  &  -1.28\,E\,+3 &  -1.23\,E\,+3          &  \textbf{-1.90\,E\,+3} \\ \cline{2-6}
		                                 & 10 &  -1.48\,E\,+3  &  -1.83\,E\,+3 &  -1.79\,E\,+3          &  \textbf{-2.67\,E\,+3} \\ \hline
        \multirow{2}{*}{Salomon}         & 5  &  3.45\,E\,0  &  2.64\,E\,0 &  2.19\,E\,0          &  \textbf{6.38\,E\,-1} \\ \cline{2-6}
		                                 & 10 &  5.05\,E\,0  &  5.79\,E\,0 &  5.29\,E\,0          &  \textbf{2.40\,E\,0} \\ \hline
        \multirow{2}{*}{Brown}           & 5  &  \textbf{6.46\,E\,-4}  &  1.39\,E\,-1 &  8.29\,E\,-2          & 3.90\,E\,-3 \\ \cline{2-6}
		                                 & 10 &  \textbf{4.76\,E\,-2}  &  9.78\,E\,-1 &  9.61\,E\,-1          &  1.65\,E\,-1 \\ \hline
	\end{tabular}
	\caption{Results summary for comparative tests: \\ averages of  $z^{*\iter{n}}$ over 100 runs, after 500 iterations per trial.}        
	\label{table:performance-results}
\end{center}
\end{table*}

\begin{table*}[!t]
    \footnotesize
    \begin{center}
    	\begin{tabular}{|c|c||c|c||c|c||c|c||c|c|}
    		\hline
    		\multirow{2}{*}{\textbf{Test function}} & \multirow{2}{*}{$D$} & \multicolumn{2}{c||}{\textbf{DIRECT}} & \multicolumn{2}{c||}{\textbf{AdaLIPO}} & \multicolumn{2}{c||}{\textbf{SMGO}} & \multicolumn{2}{c|}{\textbf{BayesianO}} \\ \cline{3-10}
    		                                 &    & 1st & 2nd & 1st & 2nd & 1st & 2nd & 1st & 2nd \\ \hline
    		                  Rosenbrock     & 10 & 3  & 34 & 0   & 0 & \textbf{75}  & 12 & 22 & \textbf{54} \\ \hline
    		\multirow{2}{*}{Styblinski-Tang} & 5  & 12 & \textbf{88} & 0   & 0 & 0   & 0 & \textbf{88} & 12 \\ \cline{2-10}
    		                                 & 10 & 15 & \textbf{71} & 0   & 0 & 8   & 15 & \textbf{77} & 14 \\ \hline
    		\multirow{2}{*}{Deb's \#1}       & 5  & 0  & 0 & 2   & \textbf{61} & 2   & 35 & \textbf{96} & 4 \\ \cline{2-10}
    		                                 & 10 & 0  & 0 & 12  & \textbf{37} & 37  & 25 & \textbf{51} & 38 \\ \hline
    		\multirow{2}{*}{Deb's \#2}       & 5  & 0  & 0 & 27  & 28 & 28  & \textbf{36} & \textbf{45} & 36 \\ \cline{2-10}
    		                                 & 10 & 0  & 0 & 29  & 28 & \textbf{44}  & 25 & 27 & \textbf{47} \\ \hline
    		\multirow{2}{*}{Schwefel}        & 5  & 0  & \textbf{82} & 1   & 8 & 0   & 9 & \textbf{99} & 1 \\ \cline{2-10}
    		                                 & 10 & 0  & 0 & 0   & \textbf{56} & 0   & 44 & \textbf{100}& 0 \\ \hline
            \multirow{2}{*}{Salomon}         & 5  & 0  & 0 & 0   & 30 & 0   & \textbf{70} & \textbf{100} & 0 \\ \cline{2-10}
    		                                 & 10 & 0  & \textbf{48} & 0   & 8 & 0   & 44 & \textbf{100}& 0 \\ \hline
            \multirow{2}{*}{Brown}           & 5  & \textbf{98}  & 2 & 0   & 0 & 0   & 0 & 2 & \textbf{98} \\ \cline{2-10}
    		                                 & 10 & \textbf{100}  & 0 & 0   & 0 & 0   & 0 & 0 & \textbf{100} \\ \hline
    	\end{tabular}
    	
    	\caption{Results summary for comparative tests: \\ number of best and 2nd best trials for the tested algorithms}
    	\label{table:best-results}
    \end{center}
\end{table*}

\begin{table*}[!t]
\footnotesize
\begin{center}
	\begin{tabular}{|c|c|c|c|c|}
		\hline
		\textbf{Test function} & $D$ & \textbf{AdaLIPO} & \textbf{SMGO} & \textbf{BayesianO} \\ \hline
		                  Rosenbrock     & 10 & 0.017 & 74.96 & 1071.04 \\ \hline
		\multirow{2}{*}{Styblinski-Tang} & 5  & 0.013 & 7.71 & 511.66 \\ \cline{2-5}
		                                 & 10 & 0.014 & 75.01 & 1085.11 \\ \hline
		\multirow{2}{*}{Deb's \#1}       & 5  & 0.012 & 7.71 & 574.08 \\ \cline{2-5}
		                                 & 10 & 0.014 & 74.64 & 1164.80 \\ \hline
		\multirow{2}{*}{Deb's \#2}       & 5  & 0.012 & 7.72 & 576.63 \\ \cline{2-5}
		                                 & 10 & 0.014 & 73.97 & 1176.10 \\ \hline
		\multirow{2}{*}{Schwefel}        & 5  & 0.011 & 7.77 & 583.24 \\ \cline{2-5}
		                                 & 10 & 0.014 & 75.13 & 1188.08 \\ \hline
		\multirow{2}{*}{Salomon}         & 5  & 0.013 & 7.63 & 643.25 \\ \cline{2-5}
		                                 & 10 & 0.015 & 75.08 & 1033.54 \\ \hline
		\multirow{2}{*}{Brown}           & 5  & 0.012 & 7.82 & 476.40 \\ \cline{2-5}
		                                 & 10 & 0.013 & 75.47 & 886.81 \\ \hline
	\end{tabular}
	\caption{Comparison of total optimization run times, averaged over the \\100 independent runs for each test case (all in seconds)}        
    \vspace{-0.25in}
	\label{table:iteration-times-comparison}
\end{center}
\end{table*}

\subsection{Computational time}
\begin{figure}[!t]
	\centering
 	\includegraphics[width=0.85\columnwidth]{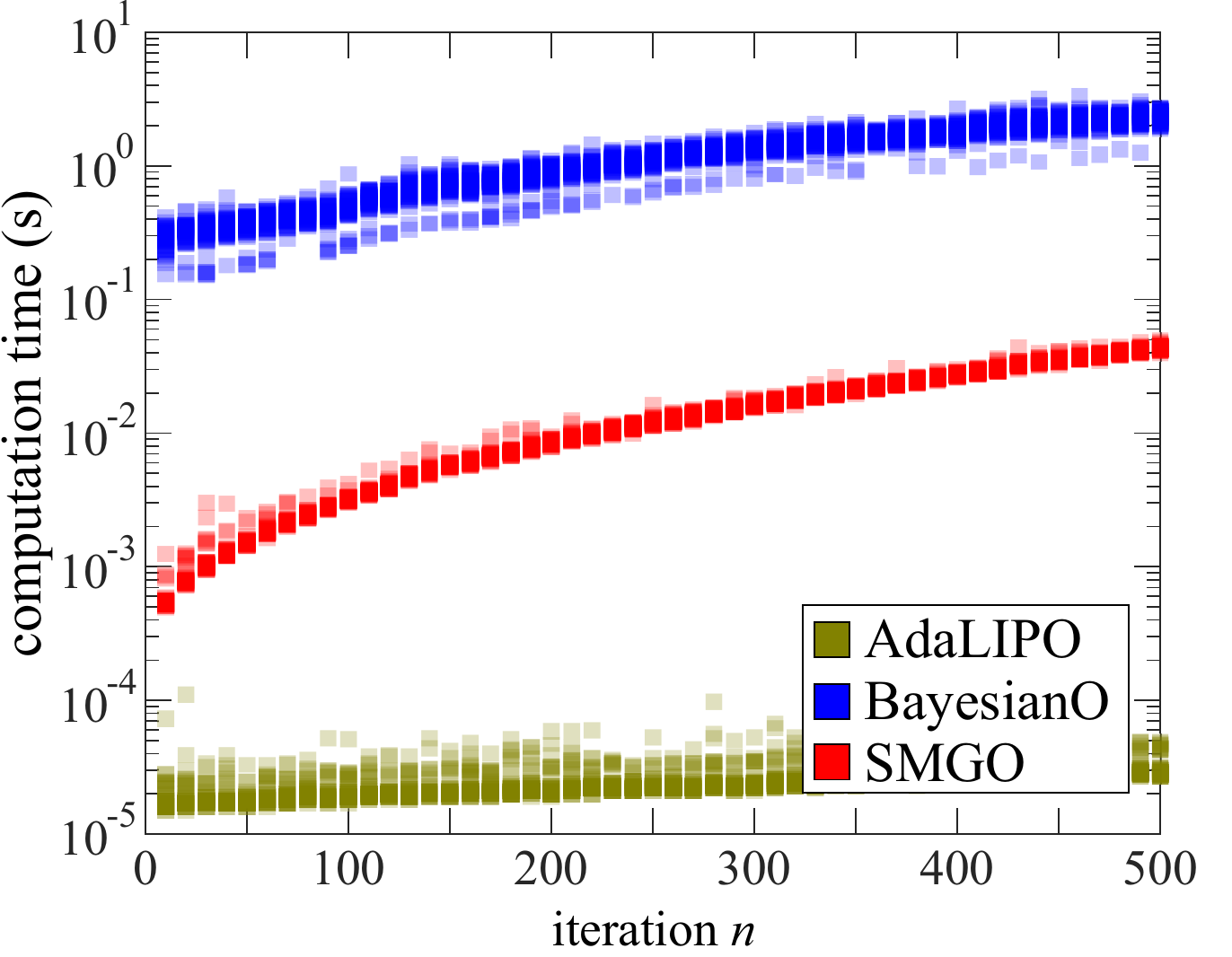}
	\caption{Calculation times on 5D Deb's \#1 (log scale)}
	\label{fig:results-calc-time}
\end{figure}

The computational time required to complete an iteration for each considered algorithms is analyzed for the case of the  5D Deb's function \#1 optimization. Fig.~\ref{fig:results-calc-time} (in log scale) shows the time required to complete an iteration as $n$ increases. The mechanism of DIRECT calculates recursive hyperboxes to be sampled by batches, hence the computational time \textit{per iteration} is not well defined, therefore  results are not shown in the plot. For AdaLIPO an almost flat line is observed, with 10 to 100~$\mu$s per iteration, implying a computational time independent of the iteration number. SMGO has resulted in a polynomial complexity w.r.t. iterations, which is expected from the computational analysis presented in Section \ref{sec:implement-notes}. Finally, BayesianO is the most demanding algorithm, requiring computational times per iteration between 2 and 3 orders of magnitude larger than SMGO in this test case.

The computational times observed for all the considered test functions  have shown  trends similar to those on the 5D Debs function \#1. Table~\ref{table:iteration-times-comparison} summarizes the results, showing the total computational time required to complete an optimization run. The AdaLIPO algorithm, based on a simple mechanism of randomized generation of exploitation and exploration points, always shows very fast computational times, not affected by the dimensionality of the decision variable, scarcely affected by $n$, and taking 0.01-0.02~s per optimization run for all the test functions.

From the tests, BayesianO requires the largest computational times, which achieved around 4-5~s per iteration at $n=500$ for most of the 10 dimensions functions. On the other hand, SMGO exhibits an intermediate per-iteration time around 100~$\mu$s to 10~ms for 5 dimensions functions and tens to hundreds of ms for 10 dimensions functions. Furthermore, for SMGO the total optimization run time when applied to 10 dimensions functions is around 10 times those employed for 5 dimensions functions. However, this is significantly faster than BayesianO, ranging from around 15 times shorter time for 10D functions, to 70 times shorter time for 5D functions.

The observed computational times, combined with the competitive optimization results presented in Figs.~\ref{fig:results-rosen-10d}-\ref{fig:results-deb1-10d}, and Table~\ref{table:performance-results}, indicate that the proposed SMGO algorithm shows a good trade-off between optimization performance and computational speed, compared to the other considered methods.

\section{Conclusions and Further Work}
\label{sec:conclusion}
In this work a sequential algorithm for global optimization of black-box functions has been proposed. The SMGO algorithm assumes a Lipschitz-continuous cost function and is based on a nonlinear Set Membership function approximation. The selection of the test points where the black-box function is evaluated is carried out by solving simplified geometry-based problems on the guaranteed lower bounds or uncertainty intervals of the function, provided by a Set Membership model. It is shown that the SMGO algorithm finds a sub-optimal solution with desired tolerance after finite iterations and converges asymptotically to the global minimum. Furthermore, the computational complexity of the algorithm is discussed and a cache-based solution is introduced to improve the computational burden.

The proposed method is evaluated against other representative global optimization methods on several test functions, comparing the quality of the best solutions found after a fixed number of calls to the black-box function. The results show the competitiveness of the SMGO algorithm compared to state of the art methods found in literature. Ongoing research aims to extend the algorithm to include constraints and its evaluation on experimental applications.

\bibliographystyle{IEEEtran}
\bibliography{smgo-paper}

\end{document}